\documentclass{gtart}

\def\ifplaintex{\expandafter\ifx\csname documentclass\endcsname\relax}

\expandafter\ifx\csname epsfbox\endcsname\relax\input epsf\fi

\ifplaintex 
\else
\fi

\def\gt{{\mathsurround=0pt\it $\cal G\mskip-2mu$eometry \&\ 
$\cal T\!\!$opology}}        

\def\gtp{{\mathsurround=0pt\it $\cal G\mskip-2mu$eometry \&\ 
$\cal T\!\!$opology $\cal P\!$ublications}}  

\def\lognumber#1{\def\thelognumber{#1}}
\def\volumenumber#1{\def\thevolumenumber{#1}}
\def\papernumber#1{\def\thepapernumber{#1}}
\def\volumeyear#1{\def\thevolumeyear{#1}}

\def\pagenumbers#1#2{\def\startpage{#1}\def\finishpage{#2}}
\def\published#1{\def\publishdate{#1}}
\def\proposed#1{\def\theproposer{#1}}
\def\seconded#1{\def\theseconders{#1}}
\def\received#1{\def\receiveddate{#1}}
\def\revised#1{\def\reviseddate{#1}}
\def\accepted#1{\def\accepteddate{#1}}
\def\asciititle#1{\def\theasciititle{#1}}

\long\def\asciiabstract#1{\long\def\theasciiabstract{#1}}
\def\asciikeywords#1{\def\theasciikeywords{#1}}

\let\\\par\let\thevolumenumber\relax\let\thepapernumber\relax
\let\thevolumeyear\relax\let\thesamplenumber\relax\let\startpage\relax
\let\finishpage\relax\let\publishdate\relax\let\receiveddate\relax
\let\reviseddate\relax\let\accepteddate\relax\let\theasciititle\relax
\let\theasciiauthors\relax
\let\theasciiabstract\relax\let\theasciikeywords\relax
\let\theasciiemail\relax\let\theshortauthors\relax\let\theshorttitle\relax

\long\def\maketitlep{   

\count0=\startpage

\gt\hfill      
\hbox to 77pt{\vbox to 0pt{\vglue -15pt\epsfbox{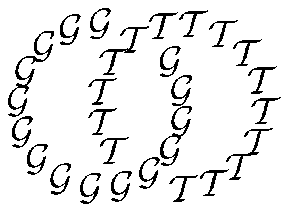}\vss}\hss}
\break
{\small\ifx\thesamplenumber\relax 
Volume \else Sample
\fi\thevolumenumber\ (\thevolumeyear)
\startpage--\finishpage\nl
Published: \publishdate}
\vglue 0.5truein plus 0.4fil minus 0.1truein

{\parskip=0pt\leftskip 0pt plus 1fil\def\\{\par\smallskip}{\ifplaintex\large
\else\Large\fi\bf\thetitle}\par\medskip}   

\vglue 0pt plus 0.1fil 

{\parskip=0pt\leftskip 0pt plus 1fil\def\\{\par}{\sc\theauthors}
\par\medskip}

\vglue 0pt plus 0.1fil 

{\small\parskip=0pt\let\newline\\
{\leftskip 0pt plus 1fil\def\\{\par}{\sl\theaddress}\par}
\expandafter\ifx\theemail\relax    
\relax\else\vglue 5pt plus 0.02fil minus 2pt\def\\{\stdspace{\rm 
and}\stdspace} 
\cl{Email:\stdspace\tt\theemail}\fi
\ifx\theurl\relax                  
\relax\else\vglue 5pt plus 0.02fil minus 2pt\def\\{\stdspace{\rm 
and}\stdspace}
\cl{URL:\stdspace\tt\theurl}\fi\par}

\vglue 7pt plus 0.3fil minus 3pt

{\bf Abstract}
\vglue 5pt plus 0.1fil minus 2pt

\theabstract

\vglue 7pt plus 0.3fil minus 3pt

{\bf AMS Classification numbers}\quad Primary:\quad \theprimaryclass

Secondary:\quad \thesecondaryclass

\vglue 5pt plus 0.3fil minus 2pt

{\bf Keywords:}\quad \thekeywords

\vglue 10pt plus 0.5fil minus 5pt

{\small  Proposed: \theproposer\hfill Received: \receiveddate\nl
Seconded: \theseconders\hfill 
\ifx\reviseddate\relax                         
Accepted: \accepteddate                        
\else
Revised: \reviseddate                          
\fi}
\eject
}       

\let\maketitlepage\maketitlep
\let\maketitle\maketitlepage

\font\phead=cmsl9 scaled 950
\font=cmsl9 scaled 1050
\font\pnum=cmbx10 scaled 913
\font\lnum=cmbx10 
\font\pfoot=cmsl9 scaled 950
\font=cmsl9 scaled 1050
\ifplaintex
\headline{\vbox to 0pt{\vskip -4.5mm\line{\small\phead\ifnum
\count0=\startpage ISSN 1364-0380 (on line)
1465-3060 (printed) \hfill {\pnum\folio}\else\ifodd\count0\def\\{ }%
\ifx\theshorttitle\relax\thetitle\else\theshorttitle\fi\hfill{\pnum\folio}
\else\def\\{ and }{\pnum\folio}\hfill\ifx\theshortauthors\relax\theauthors
\else\theshortauthors\fi\fi\fi}\vss}}
\footline{\vbox to 0pt{\vglue 0mm\line{\small\pfoot\ifnum\count0=\startpage
\copyright\ \gtp\hfill\else
\gt, Volume \thevolumenumber\ (\thevolumeyear)\hfill\fi}\vss
}}
\else
\makeatletter
\def\@oddhead{{\small\lhead\ifnum\count0=\startpage ISSN 1364-0380 (on line)
1465-3060 (printed) \hfill {\lnum\number\count0}\else\ifodd\count0
\def\\{ }\ifx\theshorttitle\relax \thetitle \else\theshorttitle\fi\hfill
{\lnum\number\count0}\else\def\\{ and }{\lnum\number\count0}
\hfill\ifx\theshortauthors\relax 
\theauthors\else\theshortauthors\fi\fi\fi}}\def\@evenhead{@oddhead}
\def\@oddfoot{\small\lfoot\ifnum\count0=\startpage\copyright\ \gtp\hfill\else
\gt, Volume \thevolumenumber\ (\thevolumeyear)\hfill\fi}
\def\@evenfoot{@oddfoot}
\makeatother
\fi

\newwrite\gtoutfile
\long\gdef\makeheadfile{  
{\def\\{, }\def\s{ }
\immediate\openout\gtoutfile head.xxx
\immediate\write\gtoutfile{To: math@arxiv.org}
\immediate\write\gtoutfile{Subject: put OR rep NNNNN:pppp}
\immediate\write\gtoutfile{--text follows this line--}
\immediate\write\gtoutfile{Proxy-for: \ifx\theasciiauthors\relax
\theauthors\else\theasciiauthors\fi\s<\ifx\theasciiemail\relax\theemail\else\theasciiemail\fi>}
\immediate\write\gtoutfile{\noexpand\\}
\immediate\write\gtoutfile{Authors: \ifx\theasciiauthors\relax
\theauthors\else\theasciiauthors\fi}
{\def\\{ }\immediate\write\gtoutfile{Title: \ifx\theasciititle\relax
\thetitle\else\theasciititle\fi}}
\immediate\write\gtoutfile{Subj-class: GT or GR or SG or ...}
\immediate\write\gtoutfile{MSC-class: \theprimaryclass\ifx\thesecondaryclass\relax\else, \thesecondaryclass\fi}
\immediate\write\gtoutfile{Journal-ref: Geom. Topol. \thevolumenumber\s
(\thevolumeyear) \startpage-\finishpage}
\immediate\write\gtoutfile{Comments: Published in Geometry and Topology at}
\immediate\write\gtoutfile{    http://www.maths.warwick.ac.uk/gt/GTVol\thevolumenumber/paper\thepapernumber.abs.html}
\immediate\write\gtoutfile{\noexpand\\}
\immediate\write\gtoutfile{}
\ifx\theasciiabstract\relax
\immediate\write\gtoutfile{\theabstract}\else
\immediate\write\gtoutfile{\theasciiabstract}\fi
\immediate\write\gtoutfile{}
\immediate\write\gtoutfile{\noexpand\\}
\immediate\write\gtoutfile{}
\immediate\closeout\gtoutfile}}  

\def\maketitlepage{\maketitlep\makeheadfile}
\let\maketitle\maketitlepage

\lognumber{170}

\volumenumber{6}
\papernumber{21}
\volumeyear{2002}
\pagenumbers{609}{647}
\received{23 March 2001}
\revised{15 March 2002}
\accepted{15 November 2002}
\proposed{Cameron Gordon}
\seconded{Walter Neumann, Michael Freedman}
\published{6 December 2002}

\usepackage{amssymb, amsmath}  
\usepackage{amscd, rlepsf}

\nocolon
\let\tilde\widetilde
\let\Tilde\widetilde

\theoremstyle{plain}
\newtheorem{theorem}{Theorem}[section]
\newtheorem{thm}{Theorem}
\newtheorem*{theorem1}{Theorem \ref{T1}}
\newtheorem*{theorem2}{Theorem \ref{T2}}
\newtheorem*{theorem3}{Theorem \ref{T3}}
\newtheorem{lemma}[theorem]{Lemma}
\newtheorem{corollary}[theorem]{Corollary}
\newtheorem{proposition}[theorem]{Proposition}

\theoremstyle{definition}
\newtheorem{definition}[theorem]{Definition}
\newtheorem{observation}[theorem]{Observation}

\theoremstyle{remark}
\newtheorem{remarks}[theorem]{Remarks}

\newtheorem*{claim*}{Claim}
\newtheorem{claim}{Claim}

\numberwithin{figure}{section}

\begin{document}

\title{Boundary curves of surfaces with the 4--plane\\property}
\asciititle{Boundary curves of surfaces with the 4-plane property}
\author{Tao Li}
\address{Department of Mathematics, Oklahoma State 
University\\Stillwater, OK 74078, USA}
\email{tli@math.okstate.edu}
\url{http://www.math.okstate.edu/\char'176tli}

\begin{abstract}
Let $M$ be an orientable and irreducible $3$--manifold whose boundary
is an incompressible torus.  Suppose that $M$ does not contain any
closed nonperipheral embedded incompressible surfaces.  We will show
in this paper that the immersed surfaces in $M$ with the $4$--plane
property can realize only finitely many boundary slopes.  Moreover, we
will show that only finitely many Dehn fillings of $M$ can yield
3--manifolds with nonpositive cubings.  This gives the first examples
of hyperbolic 3--manifolds that cannot admit any nonpositive cubings.
\end{abstract}

\asciiabstract{Let M be an orientable and irreducible 3-manifold
whose boundary is an incompressible torus.  Suppose that M does not
contain any closed nonperipheral embedded incompressible surfaces.  We
will show in this paper that the immersed surfaces in M with the
4-plane property can realize only finitely many boundary slopes.
Moreover, we will show that only finitely many Dehn fillings of M
can yield 3-manifolds with nonpositive cubings.  This gives the first
examples of hyperbolic 3-manifolds that cannot admit any nonpositive
cubings.}

\primaryclass{57M50} \secondaryclass{57M25, 57N10, 57M07}
\keywords{3--manifold, immersed surface, nonpositive cubing, 4--plane
property, immersed branched surface.}
\asciikeywords{3-manifold, immersed surface, nonpositive cubing, 4--plane
property, immersed branched surface.}

\maketitlepage

\section{Introduction}
A closed irreducible 3--manifold is called Haken if it contains a
two-sided incompressible surface.  Waldhausen has proved topological
rigidity for Haken $3$--manifolds \cite{W}, ie, if two Haken
3--manifolds are homotopically equivalent, then they are homeomorphic.
However, a theorem of Hatcher \cite{Ha} implies that, in a certain
sense, most 3--manifolds are not Haken.  Immersed $\pi_1$--injective
surfaces are a natural generalization of incompressible surfaces, and
conjecturally, $3$--manifolds that contain $\pi_1$--injective surfaces
have the same topological and geometric properties as Haken
$3$--manifolds.  Another related major conjecture in 3--manifold
topology is that any 3--manifold with infinite fundamental group
contains a $\pi_1$--injective surface.

Hass and Scott \cite{HS} have generalized Waldhausen's theorem by
proving topological rigidity for 3--manifolds that contain
$\pi_1$--injective surfaces with the $4$--plane and $1$--line
properties.  A surface in a 3--manifold is said to have the $n$--plane
property if its preimage in the universal cover of the 3-manifold is a
union of planes, and among any collection of $n$ planes, there is a
disjoint pair.  The $n$--plane property is a good way to measure the
combinatorial complexity of an immersed surface.  It has been shown
\cite{RS} that any immersed $\pi_1$--injective surface in a hyperbolic
3--manifold satisfies the $n$--plane property for some $n$.

In this paper, we use immersed branched surfaces to study surfaces
with the $4$--plane property.  Branched surfaces have been used
effectively in the studies of incompressible surfaces and laminations
\cite{FO,GO}.  Many results in 3--manifold topology (eg Hatcher's
theorem \cite{Ha}) are based on the theory of branched surfaces.  We
define an immersed branched surface in a $3$--manifold $M$ to be a
local embedding to $M$ from a branched surface that can be embedded in
some 3--manifold (see definition~\ref{D:ibs}).  Immersed branched
surfaces are also used in \cite{Li}.  Using lamination techniques and
immersed branched surfaces, we show:

\begin{thm}\label{T1}
Let $M$ be a closed, irreducible and non-Haken 3--manifold.  Then
there is a finite collection of immersed branched surfaces such that
any surface in $M$ with the 4--plane property is fully carried by an
immersed branched surface in this collection.
\end{thm}

This theorem generalizes a fundamental result of Floyd and Oertel
\cite{FO} in the theory of embedded branched surfaces.  One important
application of the theorem of Floyd and Oertel is the proof of a
theorem of Hatcher \cite{Ha}, which says that incompressible surfaces
in an orientable and irreducible 3--manifold with torus boundary can
realize only finitely many slopes.  A slope is the isotopy class of a
nontrivial simple closed curve in a torus.  We say that a surface in a
3--manifold with torus boundary can realize a slope $s$ if the
boundary of this surface consists of simple closed curves with slope
$s$ in the boundary torus of the 3--manifold. If an immersed surface
can realize a slope $s$, then it extends to a closed surface in the
closed manifold obtained by Dehn filling along the slope $s$.
However, Hatcher's theorem is not true for immersed $\pi_1$--injective
surfaces in general, since there are many 3--manifolds
\cite{B,O2,BC,Ma} in which $\pi_1$--injective surfaces can realize
infinitely many slopes, and in some cases, can realize every slope.
Using Theorem~\ref{T1}, we will show that surfaces with the $4$--plane
property are, in a sense, like incompressible surfaces.  Note that
many 3--manifolds satisfy the hypotheses in Theorems \ref{T2} and
\ref{T3}, such as hyperbolic punctured-torus bundles \cite{CJR, FH}
and hyperbolic 2--bridge knot complements \cite{HT}.

\begin{thm}\label{T2}
Let $M$ be an orientable and irreducible 3--manifold whose boundary is
an incompressible torus, and let $\mathcal{H}$ be the set of injective
surfaces that are embedded along their boundaries and satisfy the
4--plane property.  Suppose that $M$ does not contain any
nonperipheral closed (embedded) incompressible surfaces.  Then the
surfaces in $\mathcal{H}$ can realize only finitely many slopes.
\end{thm}  

Aitchison and Rubinstein have shown that if a 3--manifold has a
nonpositive cubing, then it contains a surface with the $4$--plane and
$1$--line properties \cite{AR}, and hence topological rigidity holds
for such 3--manifolds.  Nonpositive cubings, which were first
introduced by Gromov \cite{G}, are an important example of CAT(0)
structure.  A 3--manifold is said to admit a nonpositive cubing if it
is obtained by gluing cubes together along their square faces under
the following conditions: (1) For each edge, there are at least four
cubes sharing this edge; (2) for each vertex, in its link sphere, any
simple 1--cycle consisting of no more than three edges must consist of
exactly three edges, and must bound a triangle.  Mosher \cite{Mo} has
shown that if a 3--manifold has a nonpositive cubing, then it
satisfies the weak hyperbolization conjecture, ie, either it is
negatively curved in the sense of Gromov or its fundamental group has
a $\mathbb{Z}\oplus\mathbb{Z}$ subgroup.

Nonpositively cubed 3--manifolds have very nice topological and
geometric properties.  A natural question, then, is how large the
class of such 3--manifolds is.  Aitchison and Rubinstein have
constructed many examples of such 3--manifolds, and only trivial
examples, such as manifolds with finite fundamental groups, were known
not to admit such cubings.  At one time, some people believed that
every hyperbolic 3--manifold admits a nonpositive cubing.  In this
paper, we give the first nontrivial examples of 3--manifolds, in
particular, the first examples of hyperbolic 3--manifolds that cannot
admit any nonpositive cubings.  In fact, Theorem~\ref{T3} says that,
in a certain sense, most 3--manifolds do not have such a cubing.

\begin{thm}\label{T3}
Let $M$ be an orientable and irreducible 3--manifold whose boundary is
an incompressible torus.  Suppose that $M$ does not contain any closed
nonperipheral (embedded) incompressible surfaces.  Then only finitely
many Dehn fillings on $M$ can yield 3--manifolds that admit
nonpositive cubings.
\end{thm}

\noindent
\textbf{Acknowledgments}\qua  This paper is a part of my thesis.  I would
like to thank my advisor Dave Gabai, who introduced this subject to
me, for many very helpful conversations. I am also very grateful to
Yanglim Choi for a series of meetings about his thesis and for his
work on immersed branched surfaces.

\section{Hatcher's trick}\label{S2}

\begin{figure}[ht!]
\cl{\relabelbox\small
\epsfxsize3.6in\epsfbox{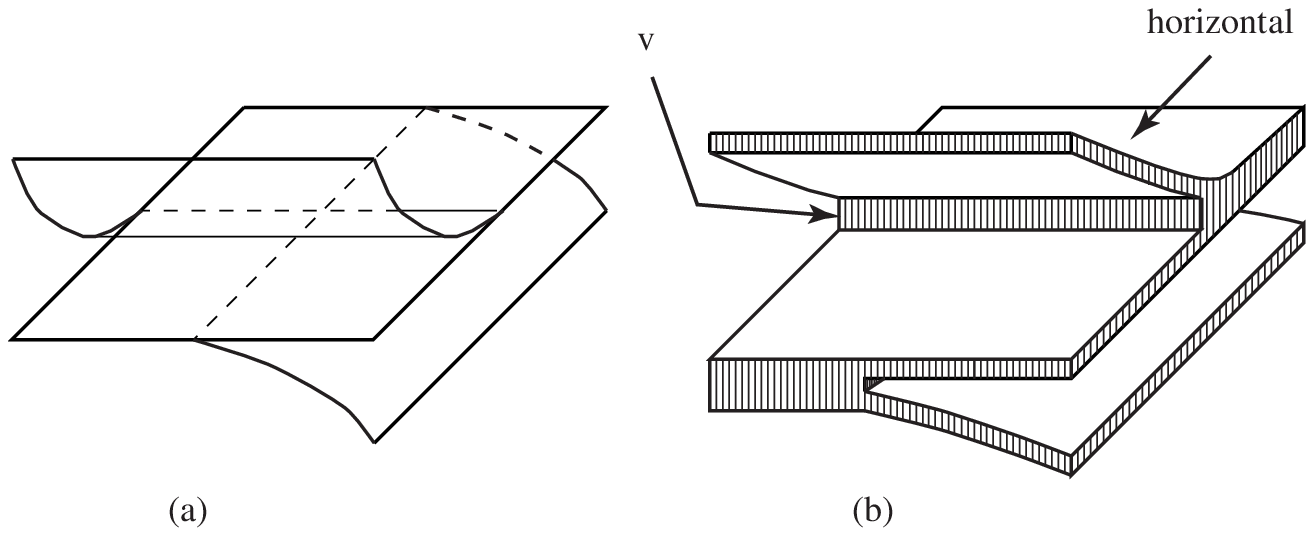}
\relabel {(a)}{(a)}
\relabel {(b)}{(b)}
\relabel {horizontal}{$\partial_hN(B)$}
\relabel {v}{$\partial_vN(B)$}
\endrelabelbox}
\nocolon
\caption{}\label{F11}
\end{figure}

A branched surface in a $3$--manifold is a closed subset locally
diffeomorphic to the model in Figure~\ref{F11}~(a).  A branched
surface is said to carry a surface (or lamination) $S$ if, after
homotopies, $S$ lies in a fibered regular neighborhood of $B$ (as
shown in Figure~\ref{F11}~(b)), which we denote by $N(B)$, and is
transverse to the interval fibers of $N(B)$.  We say that $S$ is fully
carried by a branched surface $B$ if it meets every interval fiber of
$N(B)$.  A branched surface $B$ is said to be incompressible if it
satisfies the following conditions: (1) The horizontal boundary of
$N(B)$, which we denote by $\partial_hN(B)$, is incompressible in the
complement of $N(B)$, and $\partial_hN(B)$ has no sphere component;
(2) $B$ does not contain a disk of contact; (3) there is no monogon
(see \cite{FO} for details).

\begin{theorem}[Floyd--Oertel]\label{T:FO}
Let $M$ be a compact irreducible 3--manifold with incompressible
boundary.  Then there are finitely many incompressible branched
surfaces such that every incompressible and $\partial$--incompressible
surface is fully carried by one of these branched surfaces.  Moreover,
any surface fully carried by an incompressible branched surface is
incompressible and $\partial$--incompressible.
\end{theorem}

Using this theorem and a simple trick, Hatcher has shown \cite{Ha}
that given a compact, irreducible and orientable $3$--manifold $M$
whose boundary is an incompressible torus, incompressible and
$\partial$--incompressible surfaces in $M$ can realize only finitely
many boundary slopes.  An immediate consequence of Hatcher's theorem
is that if $M$ contains no closed nonperipheral incompressible
surfaces, then all but finitely many Dehn fillings on $M$ yield
irreducible and non-Haken $3$--manifolds.  To prove Hatcher's theorem,
we need the following lemma \cite{Ha}.

\begin{lemma}[Hatcher]\label{L21}
Let $T$ be a torus and $\tau$ be a train track in $T$ that fully
carries a union of disjoint and nontrivial simple closed curves.
Suppose that $\tau$ does not bound a monogon.  Then $\tau$ is
transversely orientable.
\end{lemma}

In Theorem~\ref{T:FO}, if $\partial M$ is a torus, the boundaries of
those incompressible branched surfaces are train tracks that satisfy
the hypotheses in Lemma~\ref{L21}.  This lemma together with a trick
of Hatcher prove the following.

\begin{theorem}[Hatcher]\label{T:H}
Let $M$ be a compact, orientable and irreducible 3--manifold whose
boundary is an incompressible torus.  Suppose that $(B,\partial
B)\subset (M,\partial M)$ is an incompressible branched surface.  If
$S_1$ and $S_2$ are two embedded surfaces fully carried by $B$, then
$\partial S_1$ and $\partial S_2$ have the same slope in the torus
$\partial M$.  Moreover, the incompressible and
$\partial$--incompressible surfaces in $M$ can realize only finitely
many slopes.
\end{theorem}

\begin{figure}[ht!]
\cl{\relabelbox\small
\epsfxsize3.6in\epsfbox{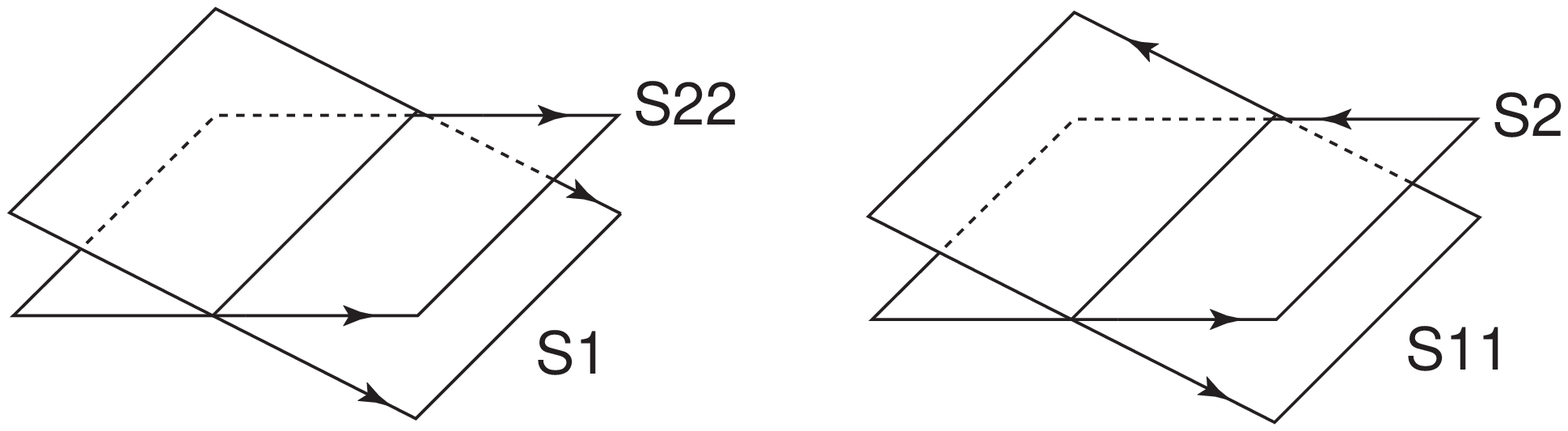}
\relabel {S1}{$S_1$} 
\relabel {S2}{$S_2$}
\relabel {S11}{$S_1$} 
\relabel {S22}{$S_2$}
\endrelabelbox}
\caption{}\label{F12}
\end{figure}

\begin{proof}
Since $M$ is orientable, the normal direction of $\partial M$ and the
transverse orientation of $\partial B$ uniquely determine an
orientation for every curve carried by $\partial B$.  Since $S_i$ is
fully carried by $B$, every component of $\partial S_i$ ($i=1\ or\ 2$)
with this induced orientation represents the same element in
$H_1(\partial M)$.  If $\partial S_1$ and $\partial S_2$ have
different slopes, they must have a nonzero intersection number.  There
are two possible configurations for the induced orientations of
$\partial S_1$ and $\partial S_2$ at endpoints of an arc $\alpha$ of
$S_1\cap S_2$, as shown in Figure~\ref{F12}.  In either case, the two
ends of $\alpha$ give points of $\partial S_1\cap\partial S_2$ with
opposite intersection numbers.  Thus, the intersection number
$\partial S_1\cdot\partial S_2=0$.  So, they must have the same slope.
The last assertion of the theorem follows from the theorem of Floyd
and Oertel.
\end{proof}

\begin{figure}[ht!]
\cl{\epsfxsize3in\epsfbox{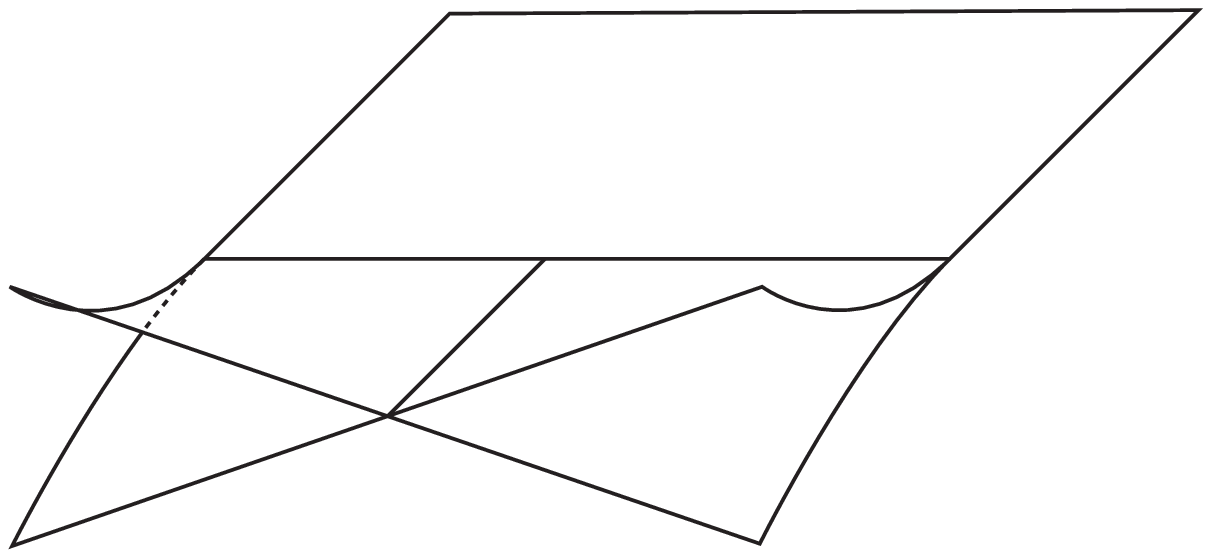}}
\caption{} \label{F13}
\end{figure}

In order to apply the trick of intersection numbers in the proof of
Hatcher's theorem, we do not need the surfaces $S_1$ and $S_2$ to be
embedded.  In fact, if $S_1$ and $S_2$ are immersed $\pi_1$--injective
surfaces that are embedded along their boundaries and transversely
intersect the interval fibers of $N(B)$, then $\partial S_1$ and
$\partial S_2$ must have the same slope by the same argument.  This is
the starting point of this paper.  In fact, even the branched surface
$B$ can be immersed.  An obstruction to applying Hatcher's trick is
the existence of a local picture as in Figure~\ref{F13} in $B$.  Next,
we will give our definition of immersed branched surfaces so that we
can apply Hatcher's trick to immersed surfaces.

\begin{definition}\label{D:ibs}
Let $B$ be a branched surface properly embedded in some compact
3--manifold, ie, the local picture of $B$ in this manifold is as in
Figure~\ref{F11}~(a).  Let $i\co B\to M$ (respectively $i\co N(B)\to
M$) be a map from $B$ (respectively $N(B)$) to a 3--manifold $M$.  We
call $i(B)$ an \emph{immersed branched surface} in $M$ if the map $i$
is a local embedding.  An immersed surface $j\co S\to M$ (or simply
$S$) is said to be \emph{carried} by $i(B)$ (or $B$) if, after some
homotopy in $M$, $j=i\circ h$, where $h\co S\to N(B)$ is an embedding
and $h(S)$ is transverse to the interval fibers of $N(B)$.  We say
that it is \emph{fully carried} by $i(B)$ if $h(S)$ transversely
intersects every $I$--fiber of $N(B)$.
\end{definition}

If $i\co B\to M$ is an immersed branched surface, then $i(B)$ contains
no local picture as in Figure~\ref{F13} by definition.  The following
proposition is an extension of Hatcher's theorem, and its proof is
simply an application of Hatcher's trick to immersed branched
surfaces.

\begin{proposition}\label{P22}
Let $M$ be a compact, orientable and irreducible 3--manifold whose
boundary is an incompressible torus.  Let $S_1$ and $S_2$ be immersed
$\pi_1$--injective surfaces fully carried by an immersed branched
surface $i\co B\to M$.  Suppose that $i|_{\partial B}$ is an embedding
and $i(\partial B)$ does not bound a monogon.  Then $\partial S_1$ and
$\partial S_2$ have the same slope.
\end{proposition} 

\section{Cross disks}\label{S3}

We have seen in section~\ref{S2} that Hatcher's trick can be applied
to immersed branched surfaces.  However, we also need finiteness of
the number of branched surfaces, as in the theorem of Floyd and
Oertel, to get interesting results.  This is impossible in general
because there are many examples of 3--manifolds in which immersed
$\pi_1$--injective surfaces can realize infinitely many slopes.  In
this section, we will show that one can generalize the theorem of
Floyd and Oertel to immersed surfaces with a certain property and such
immersed surfaces can realize only finitely many slopes.

Using normal surface theory, it is very easy to get finiteness (of the
number of branched surfaces) in the case of embedded incompressible
surfaces.  For any triangulation of a 3--manifold, an incompressible
surface can be put in Kneser--Haken normal form \cite{K,H}.  There are
7 types of normal disks in a tetrahedron, 4 triangular types and 3
quadrilateral types.  By identifying all the normal disks (in the
intersection of the surface with a tetrahedron) of the same type to a
branch sector, we can naturally construct a branched surface fully
carrying this embedded normal surface, and the finiteness follows from
the compactness of the 3--manifold (see \cite{FO} for details).
However, in the case of immersed surfaces, we cannot do this, although
immersed $\pi_1$--injective surfaces can also be put in normal
form. If we simply use the construction in \cite{FO} and identify all
the normal disks (in an immersed surface) of the same type to a branch
sector, we may get a local picture like that in Figure~\ref{F13},
which makes Hatcher's argument fail.

Suppose that $S$ is a $\pi_1$--injective surface in a 3--manifold $M$
with a triangulation $\mathcal{T}$.  Using normal surface theory, we
can put $S$ in normal form.  Let $\Tilde{M}$ be the universal cover of
$M$, $\pi \co \Tilde{M}\to M$ be the covering map,
$\Tilde{S}=\pi^{-1}(S)$, and $\Tilde{\mathcal{T}}$ be the induced
triangulation of $\Tilde{M}$.  For any arc $\alpha$ in $M$ (or
$\Tilde{M}$) whose interior does not intersect the 1--skeleton
$\mathcal{T}^{(1)}$, we define the \emph{length} of $\alpha$ to be
$|int(\alpha )\cap\mathcal{T}^{(2)}|$, where $int(E)$ denotes the
interior of $E$ and $|E|$ denotes the number of connected components
of $E$.  Moreover, we define the \emph{distance} between points $x$
and $y$, $d(x,y)$, to be the minimal length of all such arcs
connecting $x$ to $y$.  In this paper, we will always assume our
curves do not intersect the 1--skeleton of the triangulation, and we
always use the distance defined above unless specified.

Let $f\co F\to M$ be an immersed surface.  We define the \emph{weight}
of $f(F)$ to be $|f^{-1}(\mathcal{T}^{(1)})|$.  A normal (immersed)
surface $f\co F\to M$ is said to have \emph{least weight} if
$|f^{-1}(\mathcal{T}^{(1)})|$ is minimal in the homotopy class of $f$.
Let $f\co (F,\partial F)\to (M,\partial M)$ ($F\ne S^2$ or $P^2$) be a
$\pi_1$--injective map, and $M_F$ be the cover of $M$ such that
$\pi_1(M_F)$ equals $f_*(\pi_1(F))$.  We will suppose that the lift of
$f$ into $M_F$ is an embedding (note that this is automatic if $f$ is
least area in the smooth or PL sense \cite{HS, JR1}).  Thus, the
preimage of $f(F)$ in $\Tilde{M}$ consists of an embedded simply
connected surface $\Pi$ which covers $F$ in $M_F$ and the translates
of $\Pi$ by $\pi_1(M)$.  We say $f$ has the \emph{$n$--plane property}
if, given any collection of $n$ translates of $\Pi$, there is always a
disjoint pair.  We say that $\Pi$ above has \emph{least weight} if
every disk in $\Pi$ has least weight among all the disks in
$\Tilde{M}$ with the same boundary.  It follows from Theorem 5 of
\cite{JR1} or Theorem 3.4 of \cite{FHS} that $f$ can be chosen so that
$\Pi$ has least weight, and hence any translate of $\Pi$ has least
weight.  By Theorem 8 of \cite{JR1} (or Theorem 6.3 of \cite{FHS}), if
there is a map $g$ in the homotopy class of $f$ having the $n$--plane
property, then we can choose $f$ so that $f$ is a normal surface with
least weight, $\Pi$ has least weight, and $f$ also has the $n$--plane
property.  Note that $F$ may be a surface with boundary and $\Pi$ may
not be a plane $\mathbb{R}^2$, but since the interior of $\Pi$ is a
plane, to simplify notation, we will call each translate of $\Pi$ a
plane in the preimage of $f(F)$ (in $\tilde{M}$) throughout this
paper.

A \emph{normal homotopy} is defined to be a smooth map $H\co F\times
[0,1]\to M$ so that for each $t\in[0,1]$, the surface $F_t$ given by
$H|_{F\times\{ t\} }$ is a normal surface.  Note that the weight of
$F_t$ is fixed in a normal homotopy.

In this paper, we will assume that our 3--manifolds are compact and
irreducible, and our immersed surfaces, when restricted to the
boundary, are embedded.  We will also assume that our injective
surfaces are normal and have least weight, and any plane in their
preimages in the universal cover of the 3--manifold also has least
weight.  To simplify notation, we will not distinguish $f\co F\to M$,
$F$ and $f(F)$ unless necessary, and we will always denote the
preimage of $f(F)$ in the universal cover $\tilde{M}$ by $\tilde{F}$
throughout this paper.

\begin{definition} 
Let $f\co F\to M$ be a $\pi_1$--injective and least weight normal
surface, and $\Pi$ be a plane in the preimage of $f(F)$ in $\tilde{M}$
as above.  Each plane in the preimage of $f(F)$ in $\Tilde{M}$ is a
translate of $\Pi$ by an element in $\pi_1(M)$.  Let $F_1$ and $F_2$
be two such planes in $\Tilde{M}$.  Suppose that $D_1$ and $D_2$ are
two embedded subsurfaces in $F_1$ and $F_2$ respectively.  We say that
$D_1$ and $D_2$ are \emph{parallel} if there is a normal homotopy
$H\co D\times I\to\Tilde{M}$ such that $H(D,0)=D_1$, $H(D, 1)=D_2$,
$H|_{D\times\{t\}}$ is an embedding for each $t\in I$, and $H$ fixes
the 2--skeleton, ie, if $H(x,y)\in\Tilde{\mathcal{T}}^{(i)}$ then
$H(x,I)\subset\Tilde{\mathcal{T}}^{(i)}$ ($i=1,2$).  We call $D_1\cup
D_2$ a \emph{cross disk} if $D_1$ and $D_2$ are parallel disks,
$F_1\ne F_2$, and $F_1\cap F_2\ne\emptyset$.  We call $D_i$ ($i=1,2$)
a \emph{component} of the cross disk $D_1\cup D_2$.  Let $H$ be the
normal homotopy above.  We call $H(p,0)\cup H(p,1)$ a \emph{pair of
points} (respectively \emph{arcs, disks}) in the cross disk, for any
point (respectively arc, disk) $p$ in $D$.  A cross disk $D_1\cup D_2$
(or the disk $D_1$) is said to have \emph{size} at least $R$ if there
exists a point $x\in D_1$ such that $length(\alpha )\ge R$ for any
normal arc $\alpha\subset D_1$ connecting $x$ to $\partial
D_1-\partial\Tilde{M}$, and we call the normal disk of $T\cap D_1$
that contains $x$ a \emph{center} of the cross disk, where $T$ is a
tetrahedron in the triangulation.  To simplify notation, we also call
$\pi (D_1\cup D_2)$ a cross disk and call the image (under the map
$\pi$) of a pair of points (respectively arcs, disks) in $D_1\cup D_2$
a pair of points (respectively arcs, disks) in the cross disk, where
$\pi \co \Tilde{M}\to M$ is the covering map.
\end{definition}

We denote by $\mathcal{F}$ the set of $\pi_1$--injective,
$\partial$--injective and least weight surfaces in $M$ whose
boundaries are embedded in $\partial M$.  Let $\mathcal{F}_R=\{
F\in\mathcal{F}$: there are no cross disks of size $R$ in
$\Tilde{F}\}$, where $\Tilde{F}$ is the preimage of $F$ in
$\Tilde{M}$.  The following lemma is due to Choi~\cite{Ch}.

\begin{lemma}\label{L31}
There is a finite collection of immersed branched surfaces such that
every surface in $\mathcal{F}_R$ is fully carried by an immersed
branched surface in this collection.
\end{lemma} 

\begin{proof}
Let $T$ be a tetrahedron in the triangulation $\mathcal{T}$ of $M$ and
$d_i\subset F\cap T$ be a normal disk ($i=1,2,3$), where $F\in
\mathcal{F}_R$.  Suppose that $\Tilde{T}$ is a lift of $T$ in
$\Tilde{M}$, $\Tilde{d_i}$ is a lift of $d_i$ in $\Tilde{T}$, and
$F_i$ is the plane in $\tilde{F}$ that contains $\Tilde{d_i}$
($i=1,2,3$), where $\Tilde{F}$ is the preimage of $F$ in $\Tilde{M}$.
We call $D_N(d_i)=\{ x\in F_i: d(x,p)\le N,\text{ where }
p\in\Tilde{d_i}\}$ a surface of radius $N$ with center $\Tilde{d_i}$.
Note that, topologically, $D_N(d_i)$ may not be a disk under this
discrete metric.

Next, we will define an equivalence relation.  We say that $d_1$ is
\emph{equivalent} to $d_2$ if $D_{kR}(d_1)$ is parallel to
$D_{kR}(d_2)$ and $F_1\cap F_2=\emptyset$ (or $F_1=F_2$), where $k$ is
fixed.  We assume that $k$ is so large that $D_{kR}(d_i)$ contains a
subdisk of size $R$ whose center is $\Tilde{d_i}$ ($i=1,2$).  Note
that, since $M$ is compact and every plane in $\tilde{F}$ has least
weight, $k$ can be chosen to be independent of the choices of
$F\in\mathcal{F}_R$ and the normal disk $d_i\subset F$, ie, $k$
depends only on $R$ and the triangulation of $M$.  Suppose that there
are three normal disks $d_1$, $d_2$ and $d_3$ in $F\cap T$ so that
$d_1$ is equivalent to $d_2$ and $d_2$ is equivalent to $d_3$.  Then
$D_{kR}(d_1)$ is parallel to $D_{kR}(d_3)$ by definition.  If $F_1\ne
F_3$ and $F_1\cap F_3\ne\emptyset$, by the assumption on $k$, there is
a cross disk of size $R$ that consists of two disks from $F_1$ and
$F_3$.  This contradicts the hypothesis that $F\in\mathcal{F}_R$.
Thus $d_1$ is equivalent to $d_3$, and the equivalence relation is
well-defined.

Since $M$ is compact, for any normal disk $d$ in $\Tilde{M}$, the
number of nonparallel (embedded) normal surfaces of radius $kR$ (with
center $d$) is bounded by a constant $C$ that depends only on $kR$ and
the triangulation of $M$.  As there are no cross disks of size $R$, if
$D_{kR}(d_1)$ is parallel to $D_{kR}(d_2)$, then $d_1$ and $d_2$ must
be equivalent.  Thus, there are at most $C$ equivalence classes among
the normal disks of $F\cap T$ with the same disk type, and hence we
can divide the disks in $F\cap T$ for each $T$ into at most $7C$
equivalence classes, since there are $7$ different types of normal
disks in a tetrahedron.  For any tetrahedron $T$, suppose there are
$C_T$ ($C_T\le 7C$) equivalence classes in $F\cap T$.  We put $C_T$
products $D_i\times I$ ($i=1,\dots ,C_T$) in $T$ such that
$D_i\times\{ t\}$ is a normal disk and the normal disks of $F\cap T$
in the same equivalence class lie in the same product $D_i\times I$.
Along $\mathcal{T}^{(2)}$, we can glue these products $D_i\times I$'s
together according to the equivalence classes, as in the construction
of embedded branched surfaces in \cite{FO}.  In fact, we can
abstractly construct a branched surface $B$ and a map $f\co N(B)\to M$
such that, for any tetrahedron $T$,
$f(\partial_vN(B))\subset\mathcal{T}^{(2)}$ and $f(N(B)-p^{-1}(L))\cap
T$ is exactly the union of the products $int(D_i)\times I$'s in $T$,
where $L$ is the branch locus of $B$, $p\co N(B)\to B$ is the map that
collapses every interval fiber of $N(B)$ to a point, and $int(D_i)$
denotes the interior of $D_i$.  By our construction, $B$ does not
contain a local picture like that in Figure~\ref{F13}, and hence it
can be embedded in some 3--manifold \cite{C}.  Since the number of
equivalence classes is bounded by a constant, there are only finitely
many such immersed branched surfaces that fully carry surfaces in
$\mathcal{F}_R$.
\end{proof}

\begin{corollary}\label{C32}
Suppose $M$ is a compact, orientable, irreducible 3--manifold whose
boundary is an incompressible torus.  Then the surfaces in
$\mathcal{F}_R$ can realize only finitely many slopes.
\end{corollary}

\begin{proof}
Suppose that $F_1, F_2\in\mathcal{F}_R$ are fully carried by the same
immersed branched surface $f\co B\to M$.  To simplify notation, we
will also denote by $f$ the corresponding map from $N(B)$ to $M$.
Since the surfaces in $\mathcal{F}_R$ are embedded along their
boundaries, after some normal homotopy if necessary, we can assume
that $f|_{\partial B}$ is an embedding.  Since the surfaces in
$\mathcal{F}_R$ are $\pi_1$--injective, the horizontal boundary of
$f(N(B))\cap\partial M$ does not contain any trivial circle component.
Because of Lemma~\ref{L31} and Proposition~\ref{P22}, we only need to
show that $f(\partial B)$ does not bound a monogon in $\partial M$.
We will show next that the existence of a monogon in $\partial M$
contradicts our assumption that our immersed surfaces have least
weight.  The proof is essentially the same as an argument in \cite{FO}
for embedded branched surfaces.

Since $f|_{\partial B}$ is an embedding, to simplify notation, we do
not distinguish $\partial B$ and $f(\partial B)$, and denote
$f(N(\partial B))$ by $N(\partial B)$, where $N(\partial B)$ is a
fibered neighborhood of the train track $\partial B$.  By our
definition of immersed branched surface, we can assume that
$F_1\subset f(N(B))$ and $f^{-1}(F_1)$ is an embedded surface fully
carried by $N(B)$.

Suppose that $D\subset\partial M$ is a monogon, ie, $\partial
D=\alpha\cup\beta$, where $\alpha$ is a vertical arc of
$\partial_vN(\partial B)$ and $\beta\subset\partial_hN(\partial B)$.
The component of $f(\partial_v N(B))$ that contains $\alpha$ is a
rectangle $E$ whose boundary consists of two vertical arcs $\alpha,
\alpha '$ in $\partial M$ and two arcs $\gamma,\gamma'$ in
$f(\partial_vN(B)\cap\partial_hN(B))$.  Since $F_1$ is fully carried
by $f\co B\to M$, after some normal homotopy, we may assume that $E$
is embedded, $\partial_hN(\partial B)\subset\partial F_1$, and
$\gamma\cup\gamma'\subset F_1$.  Then $\delta
=\beta\cup\gamma\cup\gamma'$ is an arc in $F_1$ with
$\partial\delta\subset\partial F_1\subset\partial M$, and $\delta$ can
be homotoped rel $\partial\delta$ into $\partial M$.  Since $F_1$ is
$\partial$--injective, $\delta$ must be $\partial$--parallel in
$F_1$. So, there is an arc $\delta'\subset\partial F_1$ such that
$\delta\cup\delta'$ is a closed trivial curve in $F_1$.  Suppose
$\delta\cup\delta'$ bounds a disk $\Delta$ in $F_1$, which may not be
embedded.  Moreover, $\alpha'\cup\delta'$ also bounds a disk $D'$ in
$\partial M$, since $\alpha'\cup\delta'$ forms a homotopically trivial
curve in $M$.  So, $D\cup E\cup\Delta\cup D'$ forms an immersed sphere
in $M$.  Since $\pi_2(M)$ is trivial, we can homotope the sphere
$D\cup E\cup\Delta\cup D'$ (fixing $E$) into $E$.  After this
homotopy, we get an immersed surface in the same homotopy class as
$F_1$ with less weight.  This contradicts our least weight assumption
on the surface $F_1$.

So, $\partial B$ does not bound any monogon.  By
Proposition~\ref{P22}, $\partial F_1$ and $\partial F_2$ must have the
same slope, and the corollary follows from Lemma~\ref{L31}.
\end{proof}

\section{Limits of cross disks}\label{S4}
Let $\mathcal{H}$ be the set of injective and least weight surfaces
with the 4--plane property in $M$.  If there is a number
$R\in\mathbb{R}$ such that $\mathcal{H}\subset\mathcal{F}_R$, by
Corollary~\ref{C32}, the surfaces in $\mathcal{H}$ can realize only
finitely many slopes.  Suppose no such a number $R$ exists.  Then
there must be a sequence of surfaces $F_1, F_2,\dots
,F_n,\dots\in\mathcal{H}$ such that, in the preimage of $F_i$ in
$\Tilde{M}$ (denoted by $\Tilde{F}_i$), there is a cross disk
$D_i=D_i'\cup D_i''$ of size at least $i$, where $i\in\mathbb{N}$.
Since $M$ is compact, after passing to a subsequence if necessary, we
can assume that $D_i'$ is parallel to a subdisk $\Delta_i$ of
$D_{i+1}'$ and $d(\partial\Delta_i-\partial\Tilde{M},\partial
D_{i+1}'-\partial\Tilde{M})\ge 1$, where $d(x,y)$ denotes the
distance.  We also assume that $\partial D_i'$ lies in the
2--skeleton.

\begin{proposition}\label{P41}
The intersection of $\pi (D_i)$ with any tetrahedron does not contain
two quadrilateral normal disks of different types, where
$\pi\co\tilde{M}\to M$ is the covering map.
\end{proposition}
\begin{proof}
We know that any two quadrilateral normal disks of different types
must intersect each other.  Suppose that the intersection of $\pi
(D_i)$ with a tetrahedron contains two different types of
quadrilateral normal disks.  Let $T$ be a lift of this tetrahedron in
$\Tilde{M}$.  Then, in each of the two quadrilateral disk types, there
is a pair of parallel normal disks in $\Tilde{F}_i\cap T$ that belong
to different components of a cross disk.  By the definition of cross
disk, the two planes in $\Tilde{F}_i$ that contain the two parallel
quadrilateral normal disks must intersect each other.  Hence, the two
different quadrilateral disk types give rise to $4$ planes in
$\Tilde{F}_i$ intersecting each other.  Note that, these $4$ planes
are different planes in $\Tilde{F}_i$, since each plane is embedded in
$\tilde{M}$ by our assumptions.  This contradicts the $4$--plane
property.
\end{proof}

Thus, as in \cite{FO}, we can construct an embedded branched surface
$B_i$ in $M$ such that $\pi (D_i)$ lies in $N(B_i)$ transversely
intersecting every interval fiber of $N(B_i)$.  In fact, for each
normal disk type of $\pi (D_i)\cap T$, we construct a product
$\delta\times I$, where $T$ is a tetrahedron and $\delta\times\{ t\}$
is a normal disk of this disk type ($t\in I$).  Then, by
Proposition~\ref{P41}, we can glue these products along
$\mathcal{T}^{(2)}$ naturally to get a fibered neighborhood of an
embedded branched surface $B_i$, and $\pi (D_i)$ can be isotoped into
$N(B_i)$ transversely intersecting every interval fiber of $N(B_i)$.
Note that $B_i$ may have nontrivial boundary.  After some isotopy, we
can assume that $\partial_vN(B_i)\cap\mathcal{T}^{(1)}=\emptyset$ and
$N(B_i)\cap\mathcal{T}^{(2)}$ is a union of interval fibers of
$N(B_i)$.  By the definition of cross disk, we can also assume that
every pair of points in the cross disk lies in the same $I$--fiber of
$N(B_i)$.

\begin{proposition}\label{P42}
$N(B_i)$ can be split into an $I$--bundle over a compact surface such
that, after normal homotopies, $\pi(D_i)$ lies in this $I$--bundle,
transversely intersects its $I$--fibers, and every pair of points in
the cross disk $\pi (D_i)$ lies in the same $I$--fiber of this
$I$--bundle.
\end{proposition}
\begin{proof}
By our construction above, $N(B_i)\cap\mathcal{T}^{(2)}$, when
restricted to a 2--simplex in $\mathcal{T}^{(2)}$, is a fibered
neighborhood of a union of train tracks.  Suppose that
$\partial_vN(B_i)$ transversely intersects $\mathcal{T}^{(2)}$.
First, we split $N(B_i)$ near $N(B_i)\cap\mathcal{T}^{(2)}$ to
eliminate $\partial_vN(B_i)\cap\mathcal{T}^{(2)}$.

Let $\Delta$ be a 2--simplex in $\mathcal{T}^{(2)}$, $\delta$ be a
component of $\partial_vN(B_i)\cap\Delta$ and $N(\tau)$ be the
component of $N(B_i)\cap\Delta$ that contains $\delta$.  We associate
every such component $\delta$ of $\partial_vN(B_i)\cap\Delta$ with a
direction (in $\Delta$) that is orthogonal to $\delta$ and points into
the interior of $N(B_i)\cap\Delta$.  Let $V$ be the union of the
interval fibers of $N(\tau)$ that contain some component of
$\partial_vN(B_i)\cap\Delta$.  After performing some isotopies, we can
assume that every interval fiber in $V$ contains only one component of
$\partial_vN(B_i)\cap\Delta$.  We give every interval fiber in $V$ a
direction induced from the direction of $\partial_vN(B_i)\cap\Delta$
defined above.  Now $N(\tau) -V$ is a union of rectangles with two
horizontal edges from $\partial_hN(B_i)$ and two vertical edges from
$V$ or $\mathcal{T}^{(1)}$.  Every vertical edge from $V$ has an
induced direction.

\textbf{Case 1}\qua  For any rectangle of $N(\tau) -V$, the direction of at
most one vertical edge points inwards.

In this case, there is no ambiguity about the splitting near the
rectangle.  We split $N(\tau)$ as shown in Figure~\ref{F41}, pushing a
component of $\partial_vN(B)$ across an edge of $\Delta$.  During the
splitting we may also push some double curves of $F_i$ across this
edge.  The effect of the splitting on $\pi (D_i)$ is just an isotopy.
Thus, we can assume that any pair of points in the cross disk lies in
the same interval fiber of the fibered neighborhood of the branched
surface after this splitting.

\begin{figure}[ht!]
\cl{\relabelbox\small
\epsfxsize3.6in\epsfbox{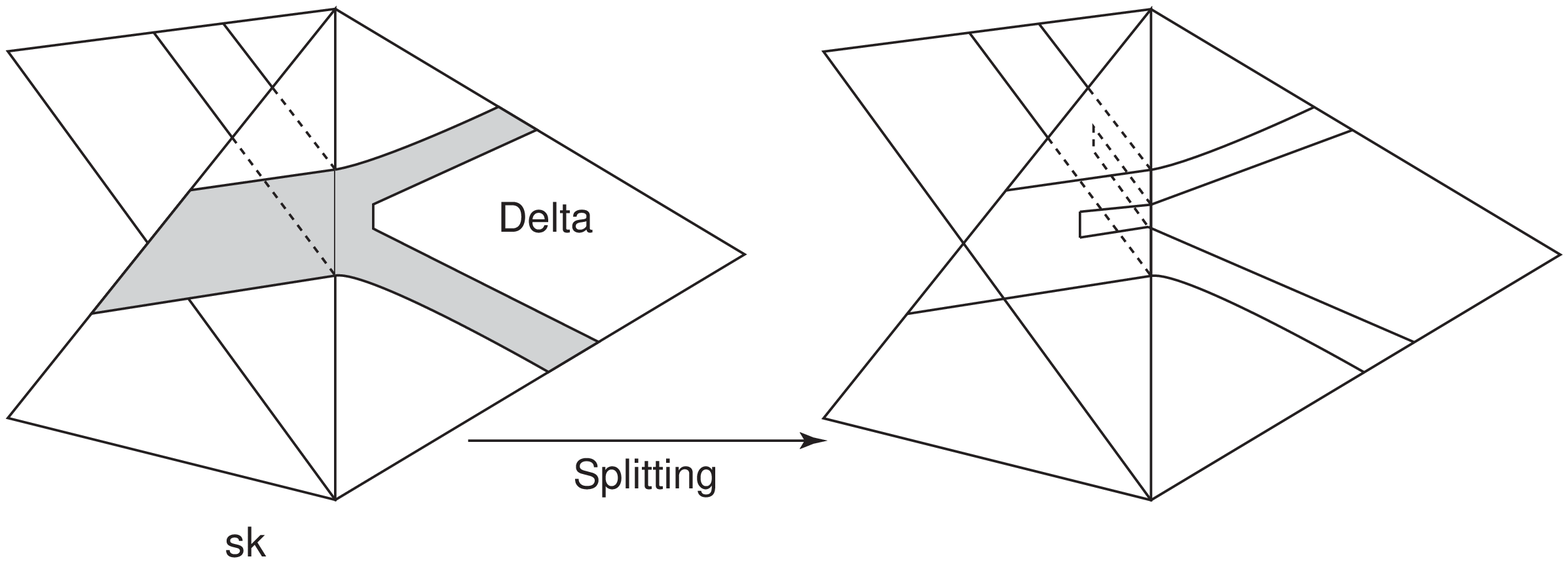}
\relabel {Delta}{$\Delta$} 
\relabel {sk}{2--skeleton}
\relabel {Splitting}{splitting} 
\endrelabelbox}
\caption{}\label{F41}
\end{figure}

\textbf{Case 2}\qua  There is a rectangle in $N(\tau) -V$ such that the
directions of both vertical edges point inwards.

The local picture of such a rectangle must be as in
Figure~\ref{F42}~(a), and there are (locally) three different
splittings as shown in Figure~\ref{F42}~(b).  We denote the rectangle
by $R$ and the part of $N(\tau)$ as in Figure~\ref{F42}~(a) by
$N(\tau)_R$.  Then $N(\tau)_R-R$ consists of 4 components, and we call
them UL (upper left) end, LL (lower left) end, UR (upper right) end
and LR (lower right) end, as shown in Figure~\ref{F42}~(a).  The
intersection of $N(\tau)_R$ and the cross disk, ie, $\pi (D_i)\cap
N(\tau)_R$, consists of arcs connecting the ends on the left side to
the ends on the right side.  An arc in $\pi (D_i)\cap N(\tau)_R$ is
called a diagonal arc if it connects an upper end to a lower end.

\begin{figure}[ht!]
\cl{\relabelbox\small
\epsfxsize3.6in\epsfbox{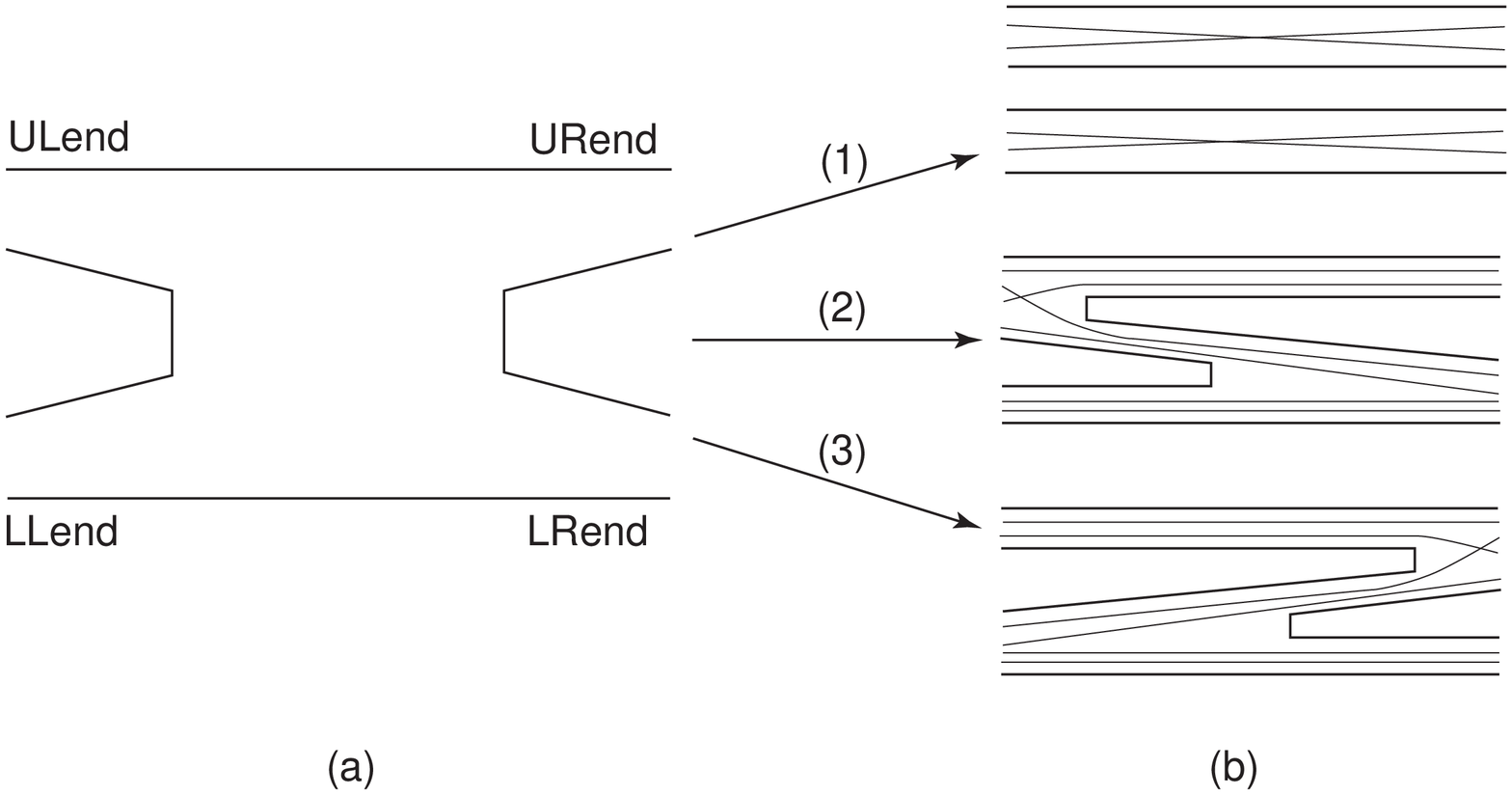}
\relabel {(a)}{(a)} 
\relabel {(b)}{(b)} 
\relabel {(1)}{(1)}
\relabel {(2)}{(2)} 
\relabel {(3)}{(3)} 
\relabel {ULend}{UL end}
\relabel {URend}{UR end} 
\relabel {LLend}{LL end} 
\relabel {LRend}{LR end}
\endrelabelbox}
\caption{}\label{F42}
\end{figure}

\begin{claim*}  
$\pi (D_i)\cap N(\tau)_R$ does not contain two diagonal arcs, say
$\alpha$ and $\beta$, such that $\alpha$ connects the UL end to the LR
end, and $\beta$ connects the LL end to the UR end.
\end{claim*}
\begin{proof}[Proof of the claim]
Suppose that it contains such arcs $\alpha$ and $\beta$.  Then there
is another arc $\alpha'$ (respectively $\beta'$) such that
$\alpha\cup\alpha'$ (respectively $\beta\cup\beta'$) is a pair of arcs
in the cross disk.  So, $\alpha'$ (respectively $\beta'$) also
connects the UL end to the LR end (respectively the LL end to the UR
end).  Note that $\alpha$ (or $\alpha'$) and $\beta$ (or $\beta'$)
must have nontrivial intersection in $N(\tau)_R$.  Next we consider a
lift of $N(\tau)_R$ in $\Tilde{M}$ and still use the same notation.
By the definition of cross disk, the 4 planes in $\Tilde{F}_i$ that
contain $\alpha$, $\alpha'$, $\beta$ and $\beta'$ respectively must
intersect each other in $\Tilde{M}$.  Since every plane in
$\Tilde{F}_i$ is embedded in $\Tilde{M}$, each is a different plane in
$\Tilde{F}_i$.  This contradicts the assumption that $F_i$ has the
4--plane property.
\end{proof}

Now we split $N(B_i)$ near $N(\tau)_R$ as follows.  If there are no
diagonal arcs in $\pi (D_i)\cap N(\tau)_R$, we split $N(B_i)$ in a
small neighborhood of $N(\tau)_R$ as the splitting (1) in
Figure~\ref{F42}.  If there are diagonal arcs, we split it as the
splitting (2) or (3) in Figure~\ref{F42} according to the type of the
diagonal arcs.  Note that by the claim, diagonal arcs of different
types cannot appear in $N(\tau)_R$ at the same time.  As in case 1, we
can assume that any pair of points of the cross disk lies in the same
$I$--fiber after the splitting.  To simplify the notation, we will
also denote the branched surface after the splitting by $B_i$.  Since
$D_i$ is compact, after finitely many such splittings,
$\partial_vN(B_i)\cap\mathcal{T}^{(2)}=\emptyset$.  Now
$\partial_vN(B_i)$ is contained in the interior of the 3--simplices,
in other words, in a collection of disjoint open 3--balls.  So, every
component of $\partial_vN(B_i)$ bounds a disk of contact.  After we
cut $N(B_i)$ along these disks of contact, as in \cite{FO},
$\partial_vN(B_i)=\emptyset$ and $N(B_i)$ becomes an $I$--bundle over
a compact surface.  As before, we can assume that, after isotopies if
necessary, every pair of points in the cross disk lies in the same
$I$--fiber.
\end{proof}

In the splittings above, we can preserve the intersection pattern of
$\Tilde{F}_i$.  For any arc $\gamma\subset F_i\cap\Delta$, since every
arc in $F_i\cap\Delta$ is a normal arc in the triangle $\Delta$, we
can assume that if an arc in $F_i\cap\Delta$ does not intersect
$\gamma$ before the splitting, it does not intersect $\gamma$ after
the splitting.  Moreover, since the intersection of $F_i$ with any
tetrahedron is a union of normal disks, we can assume that cutting the
disks of contact does not destroy the 4--plane property.  The effect
of the splitting on $F_i$ is just a normal homotopy pushing some
double curves out of the cross disk.  So, after the splitting, $F_i$
still satisfies the 4--plane property and has least weight.
Therefore, we can assume for each $i$, $\pi (D_i)$ lies in such an
$I$--bundle over a compact surface and is transverse to the
$I$--fibers.  We will still denote this $I$--bundle by $N(B_i)$.

After collapsing every $I$--fiber of $N(B_i)$ to a point, we get a
piece of embedded normal surface, which we denote by $S_i$, in $M$.
Furthermore, $D_i'$ is parallel to a subsurface of a component of
$\Tilde{S}_i$, where $\Tilde{S}_i$ is the preimage of $S_i$ in
$\Tilde{M}$.

There are only finitely many possible embedded normal surfaces (up to
normal isotopy) in $M$ that are images (under the covering map $\pi$)
of normal disks that are parallel to $D_i'$.  So, after passing to a
subsequence and doing some isotopies if necessary, we can assume that
$S_i$ is a subsurface of $S_{i+1}$.  By our assumption $d(\partial
D_i-\partial\Tilde{M}, \partial D_{i+1}-\partial\Tilde{M})\ge 1$, we
can consider the direct limit of the sequence $\{ S_i\}$ as a
(possibly noncompact) surface in $M$ whose boundary lies in $\partial
M$, and its closure is a lamination in $M$.  We can also consider this
lamination as the inverse limit of a sequence of branched surfaces
that carry $S_i$ (see \cite{M-O} for details).  We denote this
lamination by $\lambda$.  Since $\lambda$ is constructed using least
weight disks, it is well known to experts that $\lambda$ is an
essential lamination.  We provide a proof below for completeness.
Before we proceed, we will prove a useful lemma, which says that a
monogon with a long (or large) ``tail" does not exist.

\begin{lemma}\label{L41}
Let $F_0$ be a $\pi_1$--injective, $\partial$--injective least weight
normal surface in a 3--manifold $M$ and $F$ be a plane in the preimage
of $F_0$ in $\Tilde{M}$.  Suppose that $F$ has least weight and there
are two parallel disks $D_1$ and $D_2$ embedded in $F$.  Suppose that
there is a monogon, ie, an embedded disk $D$ with $\partial
D=\alpha\cup\beta$, where $\beta=D\cap (F-int(D_1\cup D_2))$,
$\beta\cap D_1$ and $\beta\cap D_2$ are the two endpoints of $\alpha$,
and $\alpha$ is an arc lying in a $2$--simplex.  Then, $weight(D_1)\le
weight(D)$.
\end{lemma}
\begin{proof}
As $D_1$ and $D_2$ are parallel, there is an embedded region
$D^2\times [1,2]$ in $\Tilde{M}$, where $D^2\times\{ t\}$ is parallel
to $D_1$ for any $t\in [1,2]$ and $D^2\times\{ i\} =D_i$ for $i=1,2$.
Moreover, by our hypothesis on $\alpha$, we can assume that $\alpha
=\{ p\}\times [1,2]$, where $p\in\partial D^2$. After some isotopy, we
can assume that $(\partial D^2\times
[1,2])\cap\tilde{\mathcal{T}}^{(1)}=\emptyset$, and hence the weight
of $\partial D^2\times [1,2]$ is zero.

We take a parallel copy of the monogon $D$, say $D'$.  Let $\partial
D'=\alpha'\cup\beta'$ and $\alpha'=\{ p'\}\times [1,2]$
($p'\in\partial D^2$), where $\alpha'$ and $\beta'$ are parallel and
close to $\alpha$ and $\beta$ respectively.  Then $\partial D^2-p\cup
p'$ consists of two arcs $\gamma$ and $\eta$.  By choosing $D'$ to be
close to $D$, we can assume that $\eta$ is the shorter one.  The four
arcs $\beta$, $\beta'$ and $\eta\times\{ 1,2\}$ form a circle that
bounds a disk $\delta$ in $F$.  We can assume that $D'$ is so close to
$D$ that the weight of $\delta$ is zero.  $D_1\cup D_2\cup\delta$ is a
disk in $F$ whose boundary is $\beta\cup\beta'\cup (\gamma\times\{
1,2\} )$. The circle $\beta\cup\beta'\cup(\gamma\times\{ 1,2\} )$ also
bounds another disk $D\cup D'\cup (\gamma\times [1,2])$ in
$\Tilde{M}$.  Since $F$ has least weight, $weight(D_1\cup
D_2\cup\delta )=2weight(D_1)\le weight(D\cup D'\cup\gamma\times
[1,2])=2weight(D)+weight(\gamma\times [1,2])$.  By our assumption
$weight(\gamma\times [1,2])=0$, we have $weight(D_1)\le weight(D)$.
\end{proof}

We call a disk as the disk $D$ in the lemma above a monogon.

\begin{lemma}\label{L42}
The lamination $\lambda$ is an essential lamination.
\end{lemma}
\begin{proof}
First we will show that every leaf of $\lambda$ is $\pi_1$--injective.
Otherwise, there is a compressing disk $D$ embedded in
$\Tilde{M}-\Tilde{\lambda }$ and $\partial D$ lies in a leaf $l$,
where $\Tilde{\lambda }$ is the preimage of $\lambda$ in the universal
cover $\Tilde{M}$.  By our construction of $\lambda$, there is, for
any $K>0$, a cross disk $D_K=D_K'\cup D_K''$ of size at least $K$ that
is parallel to a subsurface of $l$.  Since $F_K$ is $\pi_1$--injective
and has least weight, and since $\partial D$ is an essential curve in
$l$, if $K$ is large, $D_K'$ does not contain a closed curve that is
parallel to $\partial D$.  By choosing $K$ sufficiently large, we may
assume that $D_K'$ winds around $\partial D$ (in a small neighborhood
of $D$) many times, as shown in Figure~\ref{F43}~(a).  Let $N(D)$ be
an embedded disk in $\Tilde{M}$ that contains $D$ in its interior, and
$F$ be the plane in $\Tilde{F}_K$ that contains $D_K'$.  Since $F$ is
embedded in $\Tilde{M}$, the component of $F\cap N(D)$ that contains
the spiral arc in Figure~\ref{F43}~(a) must form a monogon with a long
``tail" that consists of two parallel spiral arcs winding around
$\partial D$ many times, as shown in Figure~\ref{F43}~(b).  The weight
of the monogon is at most $weight(D)$.  If $K$ is large enough, the
length of each spiral arc in the ``tail" of the monogon is very large
and, in a neighborhood of the ``tail", we can choose two parallel
disks with weight greater than $weight(D)$.  This contradicts
Lemma~\ref{L41}.

\begin{figure}[ht!]
\cl{\relabelbox\small
\epsfxsize3.6in\epsfbox{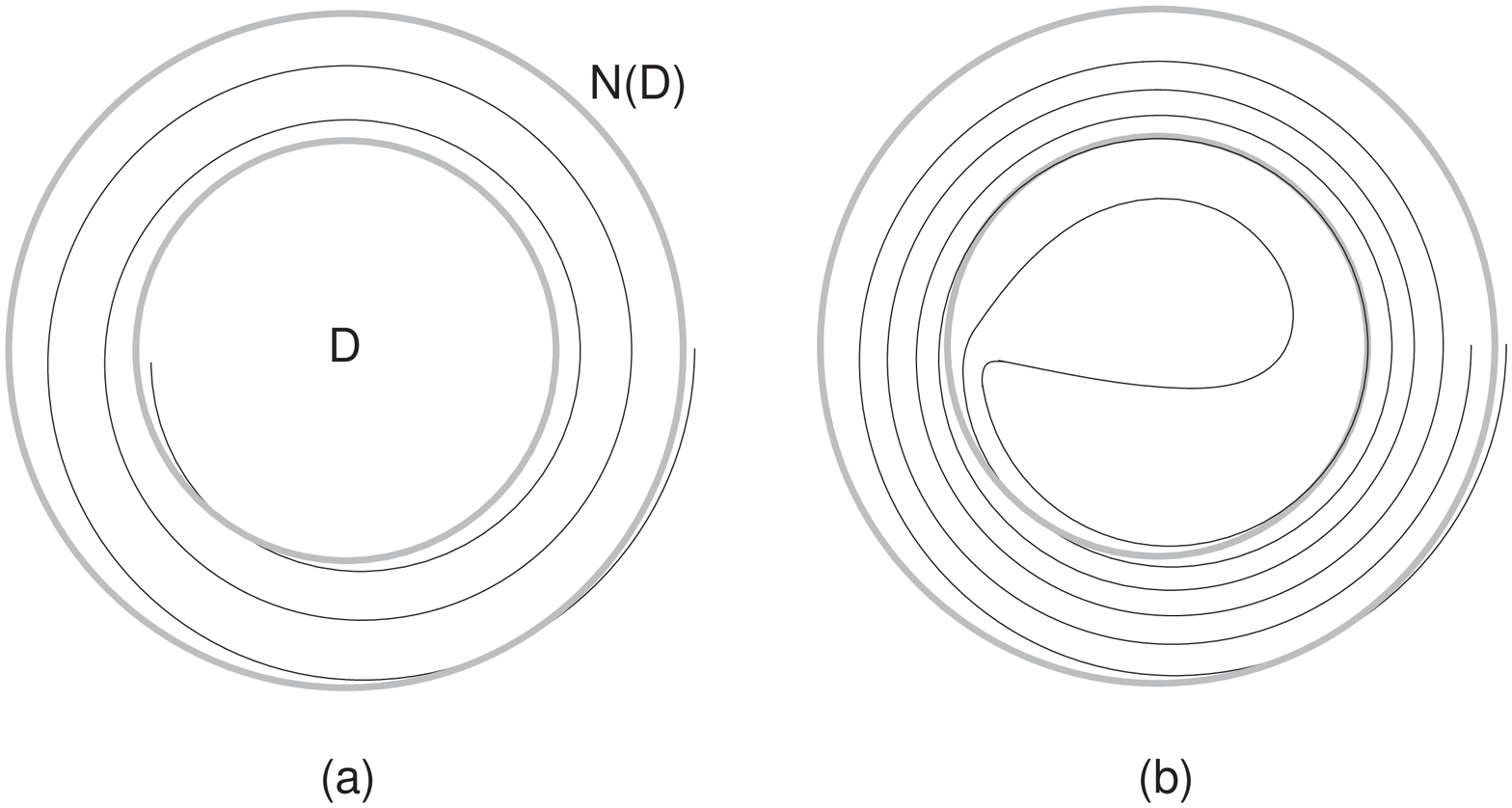}
\relabel {D}{$D$}
\relabel {N(D)}{$N(D)$} 
\relabel {(a)}{(a)}
\relabel {(b)}{(b)} 
\endrelabelbox}
\caption{}\label{F43}
\end{figure}

Next, we will show that every leaf of $\lambda$ is
$\partial$--injective.  Otherwise, there is a $\partial$--compressing
disk $D'$ whose boundary consists of two arcs $\alpha$ and $\beta$,
where $\alpha\subset\partial M$ and $\beta$ is an essential arc in a
leaf $l$.  By our construction of $\lambda$, there is a cross disk
$D_n=D_n'\cup D_n''$ of size at least $n$ such that there are arcs
$\alpha_n\subset\partial M$ and $\beta_n\subset\pi (D_n')$
($\partial\alpha_n=\partial\beta_n$) that are parallel and close to
$\alpha$ and $\beta$ respectively.  The two arcs $\alpha_n$ and
$\beta_n$ bound a disk $d_n$ that is parallel and close to $D'$.
Since the surface $F_n$ is $\partial$--injective, there must be an arc
$\gamma_n\subset\partial F_n$ such that $\gamma_n\cup\beta_n$ bounds
an immersed disk $\Delta_n$ in $F_n$.  Since $\beta$ is an essential
arc in $l$, by choosing $n$ sufficiently large, we can assume
$weight(\Delta_n)>weight(D')=weight(d_n)$.  Note that
$\gamma_n\cup\alpha_n$ must bound a disk $\delta_n$ in $\partial M$
and that $d_n\cup\Delta_n\cup\delta_n$ is an immersed $2$--sphere in
$M$.  Since $\pi_2(M)$ is trivial, we can homotope
$\Delta_n\cup\delta_n$ to $d_n$ fixing $d_n$ and get another immersed
surface $F_n'$ that is homotopic to $F_n$.  Moreover,
$weight(F_n')-weight(F_n)=weight(d_n)-weight(\Delta_n)<0$, which
contradicts the assumption that $F_n$ has least weight.

It is easy to see from our construction that no leaf is a sphere or a
disk, since the surfaces in the universal cover are embedded and are
not spheres or disks.  Also, if $\lambda$ is not end-incompressible,
there must be a monogon with a long ``tail", which contradicts
Lemma~\ref{L41} by the same argument as above.  Therefore, $\lambda$
is an essential lamination.
\end{proof}

\section{Measured sublaminations}\label{S5}

In this section, we will show that any minimal sublamination of
$\lambda$ (constructed in section~\ref{S4}) has a transverse measure.
A minimal lamination is a lamination that does not contain any proper
sublamination.  Using this result, we will prove Theorem~\ref{T1},
which can be viewed as a generalization of a theorem of Floyd and
Oertel \cite{FO}.

Let $\mu$ be a lamination in $M$ and $i\co I\times I\to M$ be an
immersion that is transverse to $\mu$, where $I=[0,1]$.  We will call
$\{ p\}\times I$ a \emph{vertical arc}, for any $p\in I$, and call
$i(I\times I)$ a \emph{transverse rectangle} if $i(I\times\partial
I)\subset\mu$ and $i^{-1}(\mu)=I\times C$ for some closed set $C$ in
$I$.

\begin{figure}[ht!]
\cl{\relabelbox\small
\epsfxsize3.6in\epsfbox{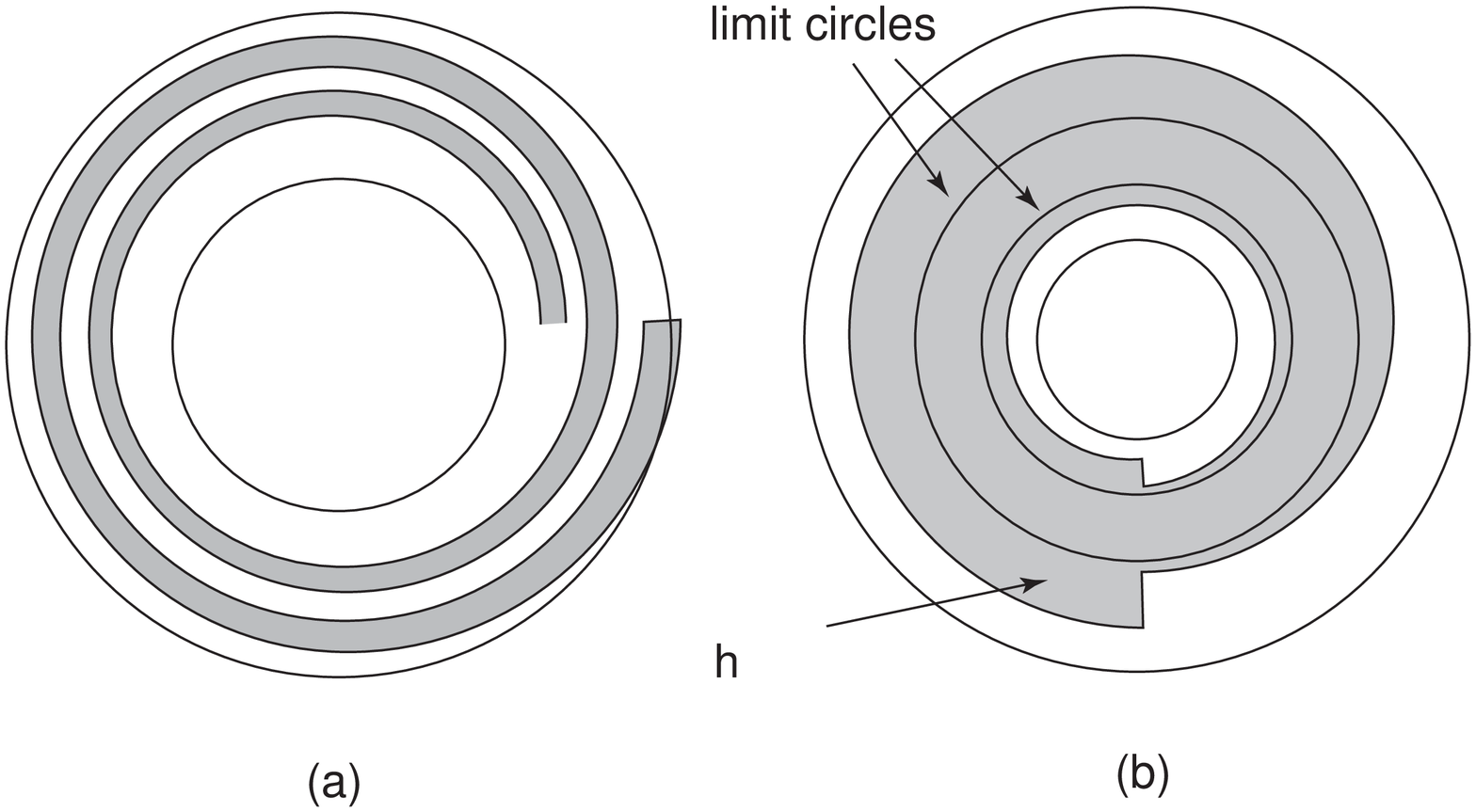}
\relabel {(a)}{(a)} 
\relabel {(b)}{(b)} 
\relabel {limit circles}{limit
circles} 
\relabel {h}{$h(I\times I)$}
\endrelabelbox}
\caption{}\label{F51}
\end{figure}

\begin{lemma}\label{L51}
Let $\mu$ be a minimal lamination.  If $\mu$ has nontrivial holonomy,
then there is a transverse rectangle $R\co I\times I\to M$ such that
$R(\{ 1\}\times I)\subset R(\{ 0\}\times int(I))$, where
$int(I)=(0,1)$.
\end{lemma}
\begin{proof}
Since $\mu$ has nontrivial holonomy, there must be a map $g\co
S^1\times I\to M$, which is transverse to $\mu$, such that
$g(S^1\times\{ 0\} )\subset L\subset\mu$ ($L$ is a leaf) and
$g^{-1}(\mu )$ consists of a collection of spirals and one circle
$S^1\times\{ 0\}$ that is the limit circle of these spirals.
Moreover, for any spiral leaf $l$ of $g^{-1}(\mu )$, there is an
embedding $i\co [0,\infty )\times I\to S^1\times I$ such that
$i^{-1}(l)=[0,\infty )\times\{ 1/2\}$ and $i([0,\infty)\times\{t\})$
is a spiral with limit circle $S^1\times\{0\}$ for each $t\in I$ (see
the shaded region in Figure~\ref{F51}~(a)).  Since $S^1\times\{ 0\}$
is the limit circle of $l$, for any arc $\{ p\}\times [0,\epsilon
]\subset S^1\times I$, there exists a number $N$, such that $i(\{
N\}\times I)\subset\{ p\}\times (0,\epsilon )$.

Since $\mu$ is a minimal lamination, every leaf is dense in $\mu$.
Thus, there is a path $\rho \co I\to L$ such that $\rho (0)=g(p,0)$,
where $p\in S^1$, and $\rho (1)\in g\circ i(\{ 0\}\times int(I))$.
Moreover, if $\epsilon$ is small enough, there is a transverse
rectangle $r\co I\times I\to M$ such that $r|_{I\times\{ 0\} }=\rho$,
$r(\{ 0\}\times I)=g(\{ p\}\times [0,\epsilon ])$, and $r(\{ 1\}\times
I)=g\circ i(\{ 0\}\times [\delta_1,\delta_2])$, where
$[\delta_1,\delta_2]\subset I$.  The concatenation of the transverse
rectangle $r$ and $g\circ i([0,N]\times [\delta_1,\delta_2])$, ie,
$R\co I\times I\to M$ where $R([0,1/2]\times I)=r(I\times I)$ and
$R([1/2,1]\times I)=g\circ i([0,N]\times [\delta_1,\delta_2])$, is a
transverse rectangle we want, where $N$ is a number that $i(\{
N\}\times I)\subset\{ p\}\times (0,\epsilon )\subset S^1\times I$.
\end{proof}

\begin{remarks}\label{R51}
$\phantom{xxx}$
\begin{enumerate}
\item The kind of construction in Lemma~\ref{L51} was also used in
\cite{Im}.
\item After connecting two copies of such transverse rectangles if
necessary, we can assume that $R(\{ 1\}\times I)\subset R(\{ 0\}\times
int(I))$ in Lemma~\ref{L51} preserves the orientation of the
$I$--fibers.  In other words, we may assume that there is a map $f\co
A\to M$ transverse to $\mu$, where $A=S^1\times I$, and an embedding
(except for $\partial I\times I$) $h\co I\times I\to A$, as shown in
Figure~\ref{F51}~(b), such that $R=f\circ h$ and $f(A)$ lies in a
small neighborhood of $R(I\times I)$.

\item Let $f$, $h$, and $R$ be the maps above. Suppose that $L_0$ and
$L_1$ are leaves in $\mu$ containing $R(I\times\{ 0\} )$ and
$R(I\times\{ 1\})$ respectively.  Then $f^{-1}(L_0\cup L_1)$ contains
two spirals of different directions whose limit circles are meridian
circles of $A$ (see Figure~\ref{F51}~(b)).  Note that $L_0$ and $L_1$
may be the same leaf and the two spirals may have the same limit
circle.

\item If $\mu$ is carried by a branched surface $B$, we can also
assume that $R(\{ q\}\times I)$ is a subarc of an interval fiber of
$N(B)$ for any $q\in I$.
\end{enumerate}
\end{remarks}

\begin{lemma}\label{L52}
Let $\lambda$ be the lamination constructed in section~\ref{S4} and
$\mu$ be any minimal sublamination of $\lambda$.  Then $\mu$ has
trivial holonomy.
\end{lemma}
\begin{proof}
Suppose that $\mu$ has nontrivial holonomy.  Since $\mu$ is a minimal
lamination, by Remarks~\ref{R51} above, there is an annulus $g\co
A=S^1\times I\to M$ such that $g^{-1}(\mu)$ contains two different
kinds of spiral leaves, as shown in Figure~\ref{F51}~(b).  From our
construction of $\lambda$, there is a cross disk $D_N=D_N'\cup D_N''$
such that $g^{-1}(\pi (D_N'))$ (respectively $g^{-1}(\pi (D_N''))$)
contains two arcs parallel and close to the two spirals respectively.
We denote these two arcs by $\alpha_0'$ and $\alpha_1'$ (respectively
$\alpha_0''$ and $\alpha_1''$), as shown in Figure~\ref{F52}~(a).  Now
we consider $g^{-1}(F_N)$, where $F_N$ is the corresponding least
weight immersed surface with the 4--plane property.  Since $F_N$ is
compact, $g^{-1}(F_N)$ is compact.  Denote the component of
$g^{-1}(F_N)$ that contains $\alpha_i'$ (respectively $\alpha_i''$) by
$c_i'$ (respectively $c_i''$), where $i=0,1$.  Since $F_N$ is a normal
surface, by Remarks~\ref{R51} (4), we can assume that $g^{-1}(F_N)$ is
transverse to each vertical arc $\{ p\}\times I$ in $A$.

\begin{figure}[ht!]
\cl{\relabelbox\small
\epsfxsize3.6in\epsfbox{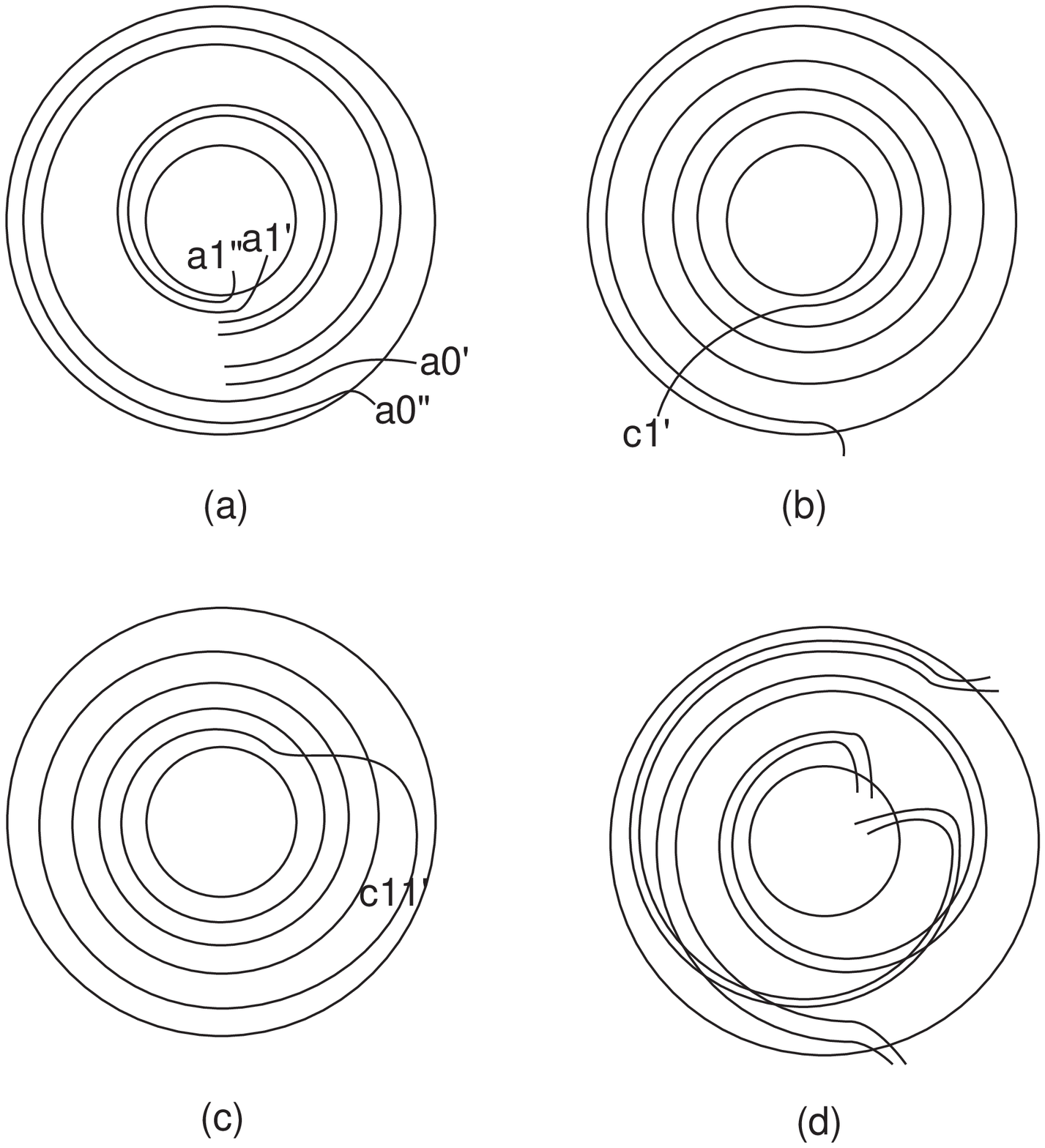}
\relabel {(a)}{(a)} 
\relabel {(b)}{(b)} 
\relabel {(c)}{(c)}
\relabel {(d)}{(d)} 
\relabel {c1'}{$c_1'$} 
\relabel {c11'}{$c_1'$} 
\relabel {a0'}{$\alpha_0'$}
\relabel {a0"}{$\alpha_0''$} 
\relabel {a1'}{$\alpha_1'$}
\relabel {a1"}{$\alpha_1''$}
\endrelabelbox}
\caption{}\label{F52}
\end{figure}

If $c_1'\cap S^1\times\{ 0\} =\emptyset$, then $c_1'$ is either a
closed curve, as shown in Figure~\ref{F52}~(c), or an arc with both
endpoints in $S^1\times\{ 1\}$, as shown in Figure~\ref{F52}~(b).
Note that, by the Reeb stability theorem, any closed curve in a leaf
with nontrivial holonomy must be an essential curve in this leaf.
Since $\lambda$ is an essential lamination, $g(S^1\times \{ 0\})$ must
be an essential curve in $M$, and we have the following commutative
diagram, where $q$ is a covering map.
$$
\begin{CD}
\mathbb{R}\times I @>\Tilde{g}>> \Tilde{M} \\ @VqVV @V\pi VV \\
A=S^1\times I @>g>> M
\end{CD}
$$
The pictures of $q^{-1}(c_1')\subset\Tilde{g}^{-1}(\Tilde{F}_N)$ are
shown in Figure~\ref{F53} (a) or (b) depending on whether $c_1'$ is an
arc with both endpoints in $S^1\times\{ 1\}$ or a closed curve.  If
$N$ is so large that $\alpha_1'$ winds around $A$ more than four
times, then there are four curves in $q^{-1}(c_1')$ intersecting each
other, as shown in Figure~\ref{F53} (a) and (b), which contradicts the
assumption that $F_N$ has the $4$--plane property.

\begin{figure}[ht!]
\cl{\relabelbox\small
\epsfxsize3in\epsfbox{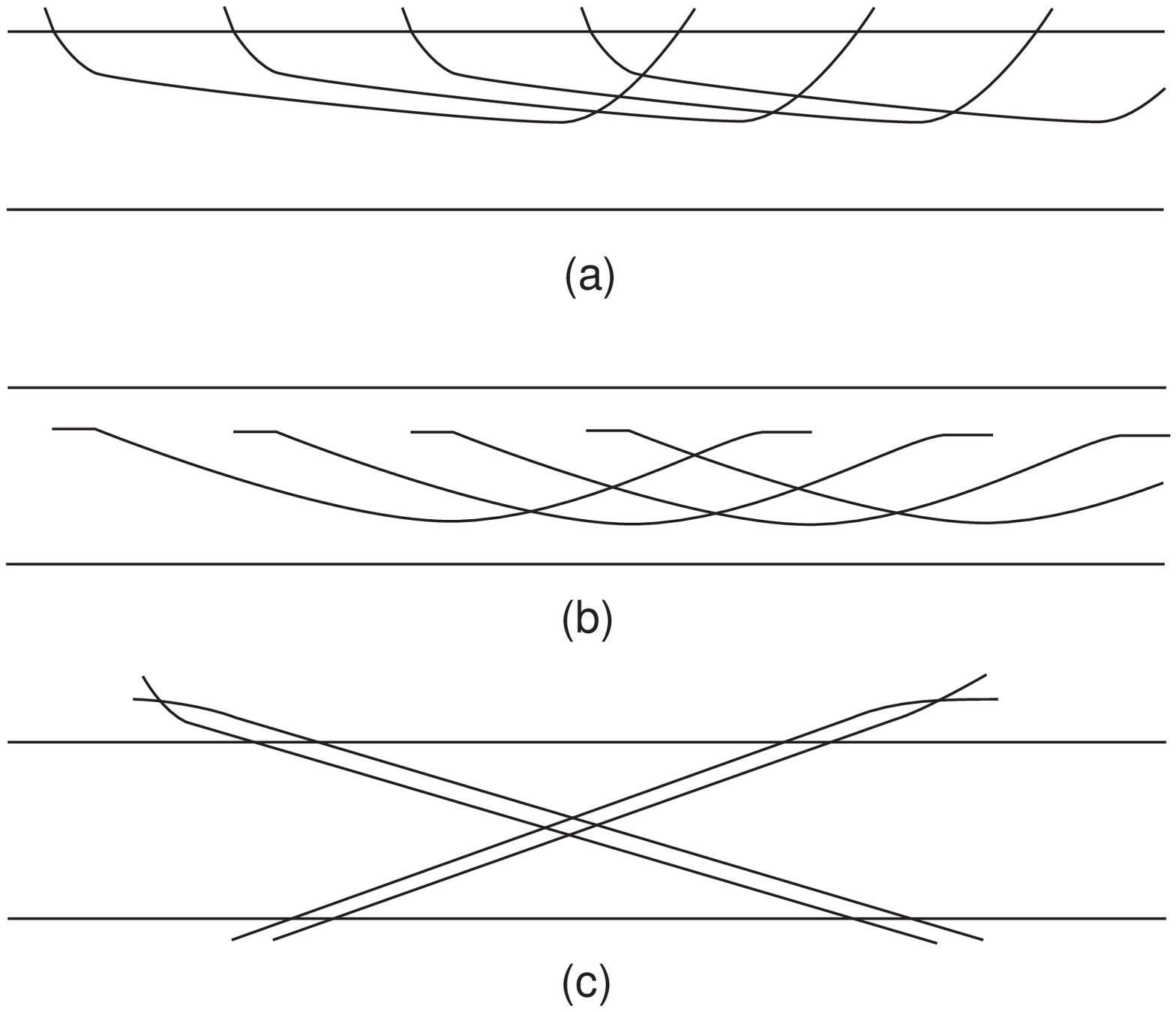}
\relabel {(a)}{(a)}
\relabel {(b)}{(b)} 
\relabel {(c)}{(c)}
\endrelabelbox}
\caption{}\label{F53}
\end{figure}

Thus, by the argument above, $c_1'$ , $c_1''$, $c_0'$ and $c_0''$ must
be arcs with endpoints in different components of $\partial A$, as
shown in Figure~\ref{F52}~(d).  In this case, $q^{-1}(c_0'\cup
c_1'\cup c_0''\cup c_1'')$ must contain $4$ arcs
$d_0',d_0'',d_1',d_1''$ as shown in Figure~\ref{F53}~(c), where
$\Tilde{g}(d_i'\cup d_i'')$ is the union of two arcs in different
components of a cross disk ($i=0,1$).  By the definition of cross
disk, the $4$ planes in $\Tilde{F}_N$ that contain $\Tilde{g}(d_0')$,
$\Tilde{g}(d_0'')$, $\Tilde{g}(d_1')$ and $\Tilde{g}(d_1'')$
respectively must intersect each other, as shown in Figure~\ref{F53}
(c), which contradicts the assumption that $F_N$ has the $4$--plane
property.
\end{proof}

The next theorem is a generalization of a theorem of Floyd and Oertel
\cite{FO}.

\begin{theorem1}
Let $M$ be a closed, irreducible and non-Haken 3--manifold.  Then
there is a finite collection of immersed branched surfaces such that
any surface in $M$ with the 4--plane property is fully carried by an
immersed branched surface in this collection.
\end{theorem1}
\begin{proof}
If the set of immersed surfaces with the $4$--plane property is a
subset of $\mathcal{F}_R$ for some number $R$ (see section~\ref{S3}
for the definition of $\mathcal{F}_R$), then the theorem follows from
by Lemma~\ref{L31}.

If there is no such a number $R$, by section~\ref{S4}, there are a
sequence of cross disks that give rise to an essential lamination
$\lambda$.  Let $\mu$ be a minimal sublamination of $\lambda$.  Since
$\mu$ is also an essential lamination, by \cite{GO}, $\mu $ is fully
carried by an embedded incompressible branched surface $B$.  By
Lemma~\ref{L52}, $\mu$ has no holonomy.  A theorem of Candel \cite{Ca}
says that if a lamination has no holonomy then it has a transverse
measure.  So, $\mu$ has a transverse measure, and hence the system of
the branch equations of $B$ (see \cite{O3}) has a positive solution.
Since each branch equation is a linear homogeneous equation with
integer coefficients, the system of branch equations of $B$ must have
a positive integer solution.  Every positive integer solution
corresponds to an embedded surface fully carried by $B$.  But, by a
theorem of Floyd and Oertel \cite{FO}, any surface fully carried by an
incompressible branched surface must be incompressible.  This
contradicts the hypothesis that $M$ is non-Haken.
\end{proof}

\section{Boundary curves}\label{S6}

Let $M$ be an irreducible 3--manifold whose boundary is an
incompressible torus, $\lambda$ be the lamination constructed in
section~\ref{S4} and $\mu$ be a minimal sublamination of $\lambda$.
Let $\{ D_i=D_i'\cup D_i''\}$ be the sequence of cross disks used in
the construction of the lamination $\lambda$ in section~\ref{S4} and
let $F_i$ be the least weight immersed surface that contains $\pi
(D_i)$.  We denote the preimage of $F_i$ in $\Tilde{M}$ by
$\Tilde{F}_i$.  Suppose that $M$ does not contain any nonperipheral
closed embedded incompressible surfaces.

\begin{lemma}\label{L61}
$\mu\cap\partial M\ne\emptyset$
\end{lemma}
\begin{proof}
Suppose that $\mu\cap\partial M=\emptyset$.  Then $\mu$ is fully
carried by an incompressible branched surface $B$ and $B\cap\partial
M=\emptyset$.  As in the proof of Theorem~\ref{T1} (see
section~\ref{S5}), the linear system of branch equations must have
integer solutions that correspond to incompressible surfaces.  Since
$B\cap\partial M=\emptyset$ and $M$ does not contain any nonperipheral
closed incompressible surfaces, those incompressible surfaces
corresponding to the integer solutions must be $\partial$--parallel
tori.

Let $N(B)$ be a fibered neighborhood of $B$, $C$ be the component of
$M-N(B)$ that contains $\partial M$, and $T_1, T_2,\dots ,T_n$ be a
collection of $\partial$--parallel tori whose union corresponds to a
positive integer solution of the system of branch equations.  After
isotopies, we can assume that every $T_i$ is transverse to the
interval fibers of $N(B)$ and $\partial_hN(B)\subset\cup_{i=1}^nT_i$.
Let $A$ be a component of $\partial_hN(B)$ that lies in the closure of
$C$.

\begin{claim*} 
The surface $A$ must be a torus.
\end{claim*}
\begin{proof}[Proof of the claim]
We first show that $A$ is not a disk.  Suppose $A$ is a disk.  Let
$\nu$ be the component of $\partial_vN(B)$ that contains $\partial A$.
Then $\partial\nu -\partial A$ is a circle in the boundary of a
component $D$ of $\partial_hN(B)$.  Since $\partial_hN(B)$ is
incompressible and $A$ is a disk, $D$ must be a disk.  So
$A\cup\nu\cup D$ is a $2$--sphere.  Since $M$ is irreducible,
$A\cup\nu\cup D$ must bound a 3--ball that contains $\cup_{i=1}^nT_i$,
which contradicts the assumption that $T_i$ is incompressible.

If $\partial A=\emptyset$, since
$\partial_hN(B)\subset\cup_{i=1}^nT_i$, $A$ is a $\partial$--parallel
torus.

Suppose $\partial A\ne\emptyset$ and $A\subset T_1$.  If there is a
component of $\partial A$ that is a trivial circle in $T_1$ then ,
since $A$ is not a disk, there must be a trivial circle in $\partial
A$ that bounds a disk in $T_1-A$.  We can isotope this disk by fixing
its boundary and pushing its interior into the interior of $N(B)$ so
that it is still transverse to the $I$--fibers of $N(B)$, and we get
a disk transverse to the $I$--fibers of $N(B)$ with its boundary in
$\partial_vN(B)$.  By definition, this is a disk of contact \cite{FO},
which contradicts the assumption that $B$ is an incompressible
branched surface.  So, every circle of $\partial A$ must be an
essential curve in $T_1$, and hence $A$ must be an annulus.

Let $c$ be a component of $\partial A$, $\nu'$ be a component of
$\partial_vN(B)$ that contains $c$, and $c'=\partial\nu'-c$ be the
other boundary component of $\nu'$.  We denote the component of
$\partial_hN(B)$ containing $c'$ by $A'$. By the argument above, $A'$
must also be an annulus.  If $A$ and $A'$ belong to different tori,
then $\nu'$ is a vertical annulus in the product region $T^2\times I$
bounded by the two tori.  This contradicts the assumptions that those
tori are $\partial$--parallel and $\partial M\subset C$.  Thus, $A$
and $A'$ must belong to the same torus $T_1$.  Then, $\nu'$ must be an
annulus in the $T^2\times I$ region bounded by $T_1$ and $\partial M$,
and $\partial\nu'\subset T_1$.  So, the vertical arcs of $\nu'$ can be
homotoped rel $\partial\nu'$ into $T_1$.  This gives rise to a monogon
and hence contradicts the assumption that $B$ is an incompressible
branched surface \cite{FO}.  Therefore, $\partial A=\emptyset$ and $A$
must be a torus.
\end{proof}

By the claim and our assumptions, $C$ must be a product region
$T^2\times I$ where $T^2\times\{ 1\} =\partial M$ and $T^2\times\{ 0\}
=A\subset\partial_hN(B)$.  Since $\mu$ is fully carried by $B$, we can
assume that $A\subset\mu$ is a leaf.  After choosing a sub cross disk
if necessary, we can assume that there is a cross disk $D_K=D_K'\cup
D_K''$ of size at least $K$ such that $\pi (D_K')$ lies in a small
neighborhood of $A$ that we denote by $T^2\times J$, where
$J=[-\epsilon ,\epsilon ]$ and $A=T^2\times\{ 0\}$.  By choosing
$\epsilon$ small enough, we can assume $T^2\times\{ t\}$ is a normal
surface for any $t\in J$.  Let $E$ be the component of $F_K\cap
(T^2\times J)$ that contains $\pi (D_K')$ and $E'$ be a component of
the preimage of $E$ in $\Tilde{M}$.  Let $F'$ be the plane in
$\Tilde{F}_K$ that contains $E'$.  So $E'$ is embedded in a region
$\mathbb{R}^2\times J$ in $\Tilde{M}$, $\partial
E'\subset\mathbb{R}^2\times\{\pm\epsilon\}$.  By choosing $\epsilon$
small enough and isotoping $F_K$, we can assume that $E'$ is
transverse to the $J$--fibers of $\mathbb{R}^2\times J$.

If $E'$ is a compact disk, then $\partial E'$ must be a circle in
$\mathbb{R}^2\times\{\pm\epsilon\}$ and $D_K$ must be in the region
bounded by $\partial E'\times J$.  So, if $K$ is large, the disk in
$\mathbb{R}^2\times\{\pm\epsilon\}$ bounded by $\partial E'$ is large.
However, if the disk bounded by $\partial E'$ is large enough, the 4
circles $g^k(\partial E')$ ($k=0,1,2,3$) must intersect each other,
where $g$ is some element in $\pi_1(\partial M)$ that acts on
$\Tilde{M}$ and fixes $\mathbb{R}^2\times J$.  This violates the
$4$--plane property, and hence $E'$ cannot be a compact disk.

Suppose that $\Tilde{F}_K\cap (\mathbb{R}^2\times\{\pm\epsilon\} )$
contains circular components.  Let $e$ be an innermost such circle and
$F_e$ be the plane in $\Tilde{F}_K$ that contains $e$.  Then $e$
bounds a disk $D$ in $\mathbb{R}^2\times\{\pm\epsilon\}$ and bounds
another disk $D'$ in $F_e$.  We can assume that
$D'\cap\pi^{-1}(T^2\times\{\pm\epsilon\} )=\partial D'$; otherwise, we
can choose $e$ to be a circle in
$D'\cap\pi^{-1}(T^2\times\{\pm\epsilon\} )$ that is innermost in $D'$.
So, $D\cup D'$ bounds a 3--ball in $\Tilde{M}$ and $\pi (D'-\partial
D')\cap (T^2\times J)=\emptyset$.  Then, we can homotope $\pi (D')$ to
$\pi (D)$ fixing $\pi (e)$.  We denote by $F_K'$ the surface after
this homotopy and denote by $F_e'$ the plane in $\Tilde{F}_K'$ (the
preimage of $F_K'$ in $\Tilde{M}$) that contains $e$.  Let $e'$ be
another component of $\pi^{-1}(\pi (e))$ and $F_{e'}$ (respectively
$F_{e'}'$) be the plane in $\Tilde{F}_K$ (respectively $\Tilde{F}_K'$)
that contains $e'$.  Since $D$ is innermost, if $F_e\cap
F_{e'}=\emptyset$, then $F_e'\cap F_{e'}'=\emptyset$.  Hence, $F_K'$
is a surface homotopic to $F_K$ and $F_K'$ also has the $4$--plane
property.  Note that since $F_K$ has least weight and $\mu$ is the
``limit" of least weight cross disks, both $D$ and $D'$ have least
weight and $weight(D)=weight(D')$.  Thus, $F_K'$ also has least weight
and $F_K'\cap T^2\times\{\pm\epsilon\}$ has fewer trivial circles
after a small homotopy.  So, we can assume that
$\Tilde{F}_K\cap\mathbb{R}^2\times\{ \pm\epsilon\}$ contains no
trivial circles.  Note that since $E'$ can never be a compact disk by
the argument above, this homotopy will not push the entire $E'$ out of
$\mathbb{R}^2\times J$.  Therefore, we can assume that $E'$ is a
noncompact and simply connected surface.

If $\partial E'\cap\mathbb{R}^2\times\{\epsilon\}$ has more than one
component, then since we have assumed that $E'$ is transverse to the
$J$--fibers of $\mathbb{R}^2\times J$, $\partial
E'\cap\mathbb{R}^2\times\{\epsilon\}$ bounds a (noncompact) region $Q$
in $\mathbb{R}^2\times\{\epsilon\}$, $D_K'\subset Q\times J$, and
$\partial Q$ contains more than one line.  Moreover, since $E'$ is
transverse to the $J$--fibers, it is easy to see that, for any element
$g\in\pi_1(\partial M)$ that acts on $\Tilde{M}$ fixing
$\mathbb{R}^2\times J$, if $Q\ne g(Q)$ and $Q\cap g(Q)\ne\emptyset$ in
$\mathbb{R}^2\times\{\epsilon\}$, then $E'\cap g(E')\ne\emptyset$.  If
$K$ is large, the distance between any two lines in $\partial Q$ must
be large.  Thus, by assuming $D_K'$ to be large, we can always find a
nontrivial element $g$ in $\pi_1(\partial M)$ such that the $g^k(Q)$'s
($k=0,1,2,3$), and hence the $g^k(E')$'s ($k=0,1,2,3$) intersect each
other, which contradicts the $4$--plane property.

Therefore, $\partial E'\cap\mathbb{R}^2\times\{\epsilon\}$ must be a
single line, and hence $E$ must be an immersed annulus in $T^2\times
J$ with one boundary component in $T^2\times\{\epsilon\}$ and the
other boundary component in $T^2\times\{ -\epsilon\}$.  By our
construction, $weight(E)$ is large if $K$ is large.  We can always
find an immersed annulus $A_E\subset T^2\times J$ with $\partial
A_E=\partial E$ and $weight(A_E)$ relatively small.  So, the surface
$(F_K-E)\cup A_E$ is homotopic to $F_K$ and has less weight.  The
homotopy is like a Dehn twist unwrapping $E$ to get $A_E$.  This
contradicts the assumption that $F_K$ has least weight in its homotopy
class.  So, $\mu\cap\partial M$ cannot be empty.
\end{proof}

\begin{lemma}\label{L62}
$\partial\mu$ is a lamination by circles.
\end{lemma}
\begin{proof}
Since $\mu$ is a measured lamination and $\partial M$ is a torus,
$\partial\mu$ is either a lamination by circles or a lamination by
lines with an irrational slope. Suppose $\mu$ is fully carried by an
incompressible branched surface $B$. Let $\mathcal{S}$ be the solution
space of the system of branch equations of $B$.  Since the
coefficients of the branch equations are integers, there are finitely
many positive integer solutions that generate $\mathcal{S}$, ie, any
point (solution) in $\mathcal{S}$ can be written as a linear
combination of these integer solutions.  Every such integer solution
gives rise to an incompressible surface fully carried by $B$.  By
Hatcher's theorem, these surfaces have the same boundary slope.  The
boundary slope of any measured lamination $\mu$ fully carried by $B$
is equal to the measure of a longitude of $\partial M$ divided by the
measure of a meridian.  Hence, the boundary slope can be expressed as
a fraction with both numerator and denominator homogeneous linear
functions of the weights of the branch sectors.  Note that, similar to
the proof of Hatcher's theorem, we can choose a transverse orientation
for $\partial B$ and assume the homogeneous linear functions above are
fixed in the calculation of the boundary slopes of any surfaces or
measured laminations fully carried by $B$.  Since the solution in
$\mathcal{S}$ that corresponds to $\mu$ is a linear combination of
those integer solutions, and since the boundary slopes of those
integer solutions (plugging into the fraction described above) are the
same, $\partial\mu$ must have the same slope as the boundary slope of
these incompressible surfaces.  Therefore, the boundary of any
measured lamination fully carried by $B$ is a lamination by circles
with the same slope.
\end{proof}

\begin{lemma}\label{L63}
Let $\{ D_i=D_i'\cup D_i''\}$ be the sequence of cross disks used in
the construction of an essential lamination in section~\ref{S4}, and
$F_i$ be the immersed surface with the $4$--plane property that
contains $\pi (D_i)$.  Then, $\{F_i\}$ contains a subsequence of
surfaces with the same boundary slope.
\end{lemma}
\begin{proof}
Let $\lambda$ be the essential lamination constructed using $\{D_i\}$
as in section~\ref{S4}, and $\mu$ be a minimal sublamination of
$\lambda$.  Then, by Lemma~\ref{L62}, $\partial\mu$ is a lamination by
circles.  Let $B$ be an incompressible branched surface that fully
carries $\mu$.  Since $\partial\mu$ is a union of parallel circles, we
can assume that $\partial B$ is a union of circles.  Let $N(B)$ be a
fibered neighborhood of $B$, $\Tilde{B}=\pi^{-1}(B)$ and
$N(\Tilde{B})=\pi^{-1}(N(B))$.  We can assume that each cross disk
$D_i$ lies in $N(\Tilde{B})$, otherwise, we can choose a large sub
cross disk of $D_i$ that lies in $N(\Tilde{B})$ for each $i$, and the
proof is the same.

Suppose the lemma is not true, then we can choose $\{F_i\}$ to be a
sequence of surfaces no two of which have the same boundary slopes.
We can also assume that $\partial F_k$ has a different slope from
$\partial\mu$ for each $k$.  Then $\pi (D_k)$ is a piece of immersed
surface in $N(B)$ transverse to every $I$--fiber, and $\pi
(D_k)\cap\partial M$ is a union of spirals in $N(B)\cap\partial M$.
We give each component of $\partial B$ an orientation so that they
represent the same element in $H_1(\partial M)$.  This orientation of
$\partial B$ determines an orientation for each $I$--fiber of
$N(B)\cap\partial M$.  As in the proof of Hatcher's theorem, the
orientation of the $I$--fibers and a normal direction of $\partial M$
uniquely determine an orientation for every curve in $N(B)\cap\partial
M$ that is transverse to the $I$--fibers of $N(B)$.

\begin{claim}\label{claim1}
If $k$ is sufficiently large, we can assume that each circle in
$\partial F_k$ admits a direction along the curve that agrees with the
induced orientation of every arc in $\partial F_k\cap N(B)$ described
above.
\end{claim}
\begin{proof}[Proof of claim \ref{claim1}]
Suppose there is a circle in $\partial F_k$ that does not admit such
an orientation.  Then there must be a subarc $C$ of the circle outside
$N(B)\cap\partial M$ connecting two spirals that are either in the
same component of $N(B)\cap\partial M$, as shown in
Figure~\ref{F61}~(a), or in different components of $N(B)\cap\partial
M$ with incompatible induced orientations, as shown in
Figure~\ref{F61}~(b).  We will show that both cases contradict our
assumption that $F_k$ is of least weight in its homotopy class.  After
assuming the size of the cross disk to be large, we can rule out the
first possibility, ie, Figure~\ref{F61}~(a), by Lemma~\ref{L41}.  To
eliminate the second possibility, ie, Figure~\ref{F61}~(b), we use a
certain triangulation of $M$ as follows.

\begin{figure}[ht!]
\cl{\relabelbox\small
\epsfxsize3.6in\epsfbox{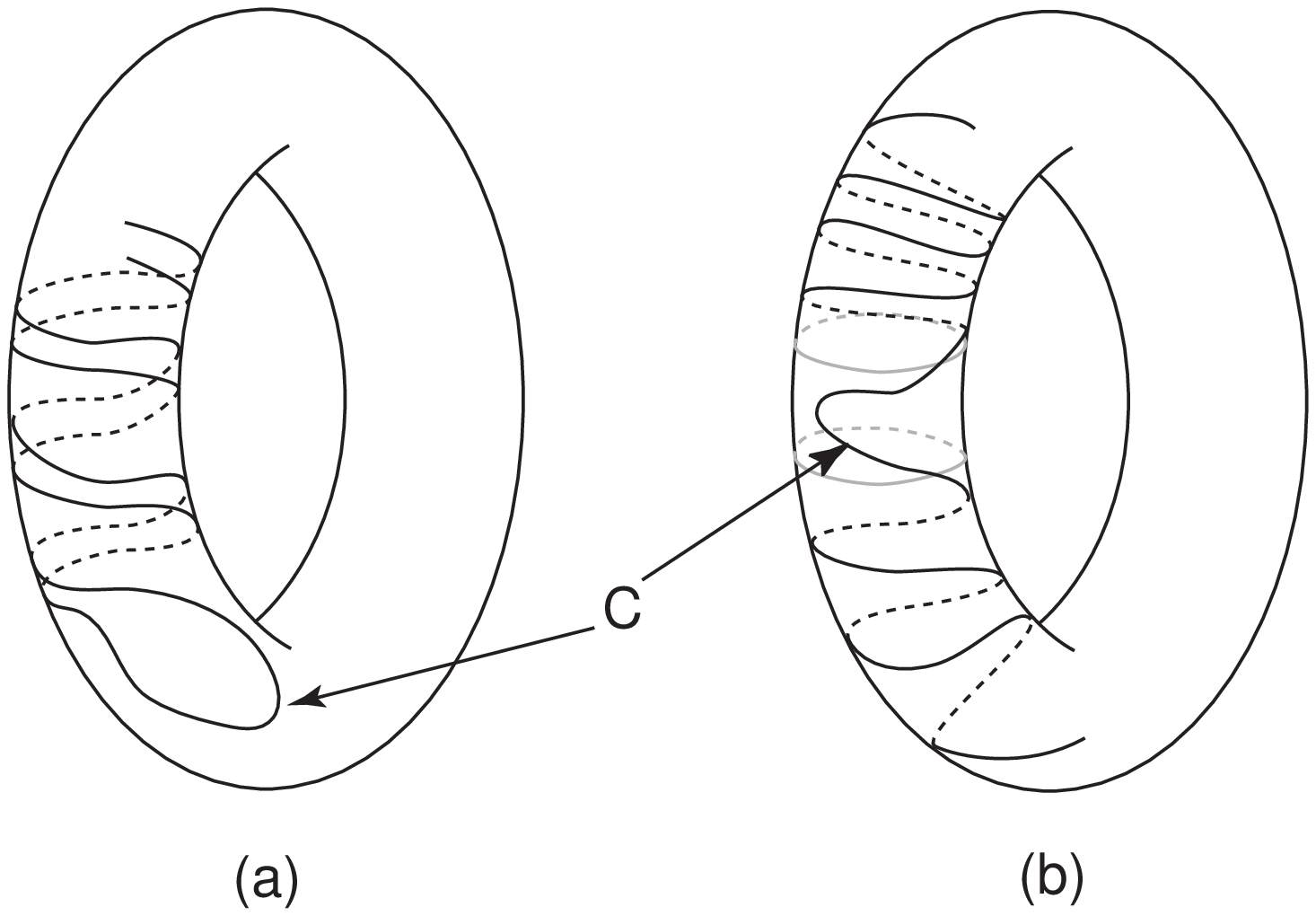}
\relabel {(a)}{(a)} 
\relabel {(b)}{(b)} 
\relabel {C}{$C$}
\endrelabelbox}
\caption{}\label{F61}
\end{figure}

By \cite{JR2}, there is a one-vertex triangulation $\mathcal{T}$ of
$M$ and this vertex is in $\partial M$.  Since $\partial M=T^2$, the
induced triangulation of $\partial M$ must consist of two triangles as
shown in Figure~\ref{F62}~(a).  Now we glue a product region
$T^2\times I$ ($I=[0,1]$) to $M$ with $T^2\times\{ 0\}=\partial M$.
Hence, $(\mathcal{T}^{(1)}\cap\partial M)\times I$ gives a cellulation
of $T^2\times I$ that consists of a pair of triangular prisms.  Then,
we add a diagonal to each rectangular face of the prisms, which gives
a triangulation of $T^2\times I$.  Figure~\ref{F62}~(b) is a picture
of the induced triangulation of a fundamental domain in the universal
cover of $T^2\times I$.  Since $M\cup (T^2\times I)$ is homeomorphic
to $M$, we can assume that $M$ has a triangulation as that of $M\cup
(T^2\times I)$ described above.  To simplify notation, we still use
$\mathcal{T}$ to denote this new triangulation of $M$.  Now,
$\mathcal{T}^{(0)}\cap\partial M$ is a single vertex $v$ and the
intersection of its link hemisphere $H$ and $\mathcal{T}^{(1)}$
consists of 10 points of which $6$ points lie in $\partial
H\subset\partial M$.

\begin{figure}[ht!]
\cl{\relabelbox\small
\epsfxsize3.6in\epsfbox{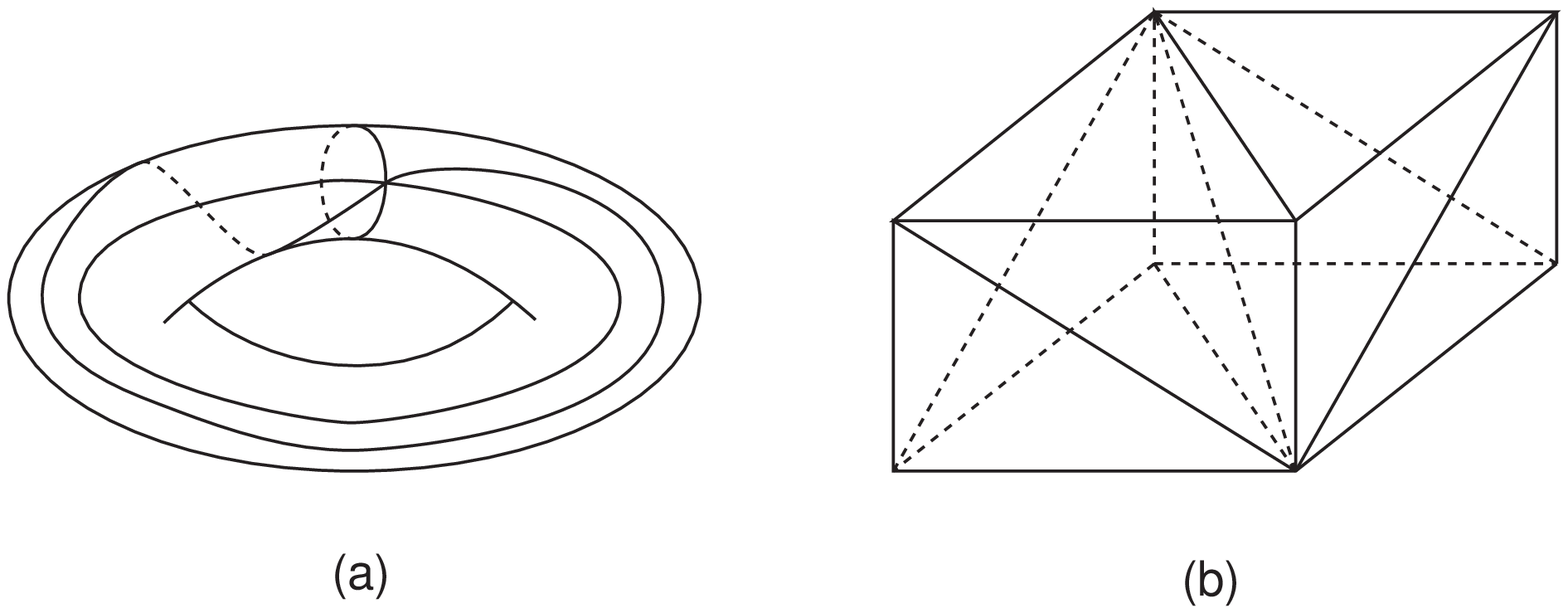}
\relabel {(a)}{(a)} 
\relabel {(b)}{(b)}
\endrelabelbox}
\caption{}\label{F62}
\end{figure}

We assume that our immersed surfaces are normal and have least weight
with respect to the triangulation above.  Suppose the second case,
ie, Figure~\ref{F61}~(b), occurs.  Let $A$ be the annular component
of $\partial M-N(B)$ that contains the arc $C$.  Then we isotope $F_k$
by pushing $C$ along $A$ to ``unwrap" the spirals in a small
neighborhood of $\partial M$, as shown in Figure~\ref{F63} (a) and
(b).  If the vertex $v$ is not in $A$, then after this isotopy,
$|\partial F_k\cap\mathcal{T}^{(1)}|$ decreases and $|(F_k-\partial
F_k)\cap\mathcal{T}^{(1)}|$ does not change.  This contradicts the
assumption that $F_k$ has least weight.  So $v\in A$.  If every edge
of $\mathcal{T}^{(1)}\cap\partial M$ intersects $\partial A$
nontrivially, then after $C$ passes through the vertex $v$ during the
isotopy, $|\partial F_k\cap\mathcal{T}^{(1)}|$ decreases by 6 and
$|(F_k-\partial F_k)\cap\mathcal{T}^{(1)}|$ increases by 4.  Hence,
the total weight of $F_k$ decreases, which also gives a contradiction.
Therefore, there is an edge $e$ of $\mathcal{T}^{(1)}\cap\partial M$
lying inside $A$, as shown in Figure~\ref{F63}~(a).  Then by our
construction of the triangulation, $e$ forms a meridian circle of the
annulus $A$ and there is at most one such edge.  After $C$ passes
through $v$ in the isotopy above, $|\partial
F_k\cap\mathcal{T}^{(1)}|$ decreases by 4, $|(F_k-\partial
F_k)\cap\mathcal{T}^{(1)}|$ increases by 4, and the total weight does
not change.

\begin{figure}[ht!]
\cl{\relabelbox\small
\epsfxsize4in\epsfbox{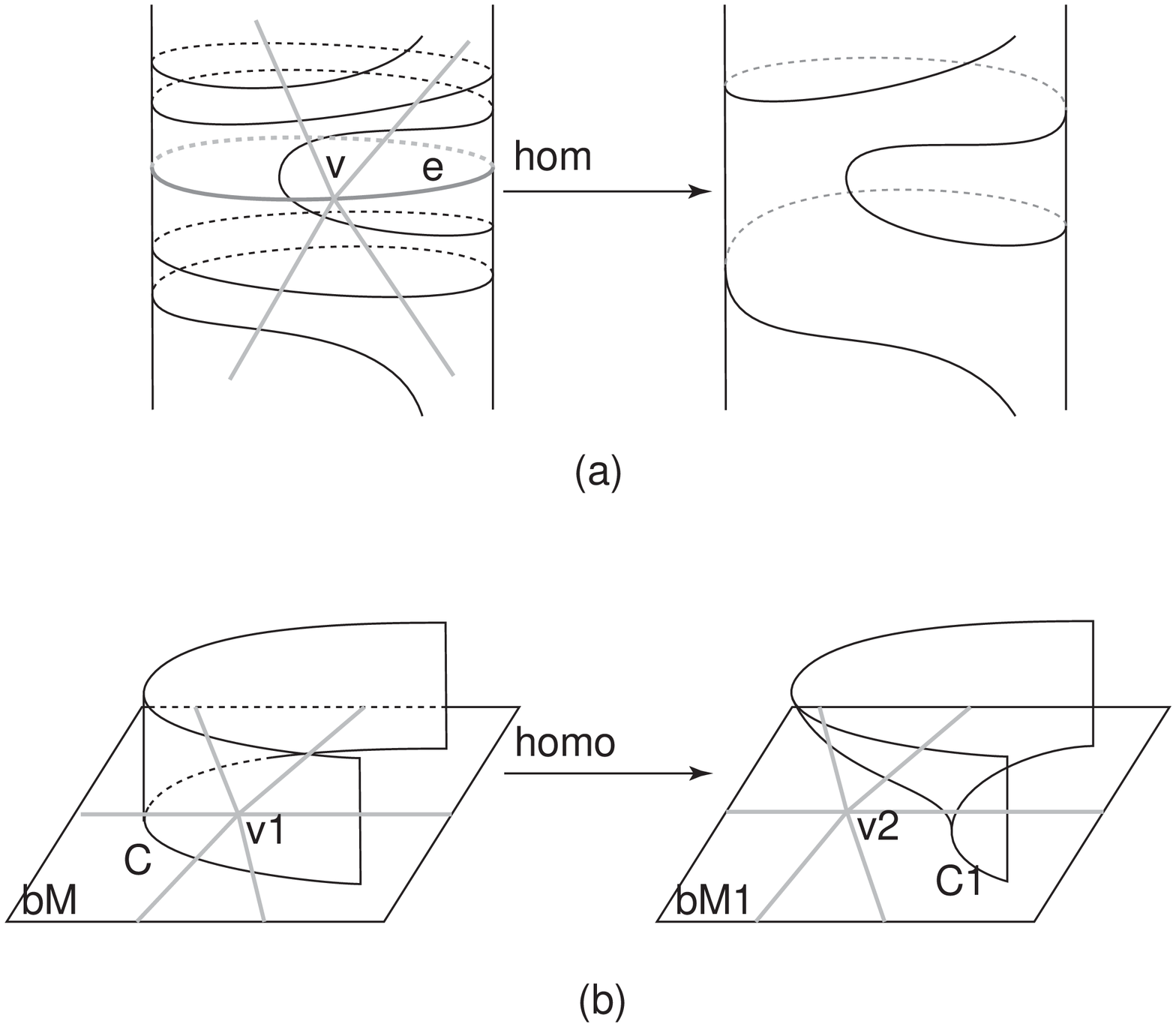}
\relabel {(a)}{(a)} 
\relabel {(b)}{(b)} 
\relabel {v}{$v$} 
\relabel {v1}{$v$} 
\relabel {v2}{$v$} 
\relabel {e}{$e$}
\relabel {hom}{homotopy} 
\relabel {homo}{homotopy} 
\relabel {C}{$C$} 
\relabel {C1}{$C$} 
\relabel {bM}{$\partial M$}
\relabel {bM1}{$\partial M$}
\endrelabelbox}
\caption{}\label{F63}
\end{figure}

Now, we will see exactly what happens in a tetrahedron.  Let $T$ be a
tetrahedron with a face $\Delta$ in $\partial M$.  There is a normal
arc $\delta$ in $C\cap\Delta$ that cuts off a subtriangle (in
$\Delta\cap A$) that contains the vertex $v$.  The normal disk of
$F_k\cap T$ containing $\delta$ is either a triangle or a
quadrilateral.  If we do the isotopy as in Figure~\ref{F63} (b) by
pushing $C$ across $v$, then the effect of this isotopy on the normal
disk that contains $\delta$ is either as in Figure~\ref{F64}~(a), in
which case the normal disk is a triangle, or as in
Figure~\ref{F64}~(b), in which case the normal disk is a
quadrilateral.  In the first case, as shown in Figure~\ref{F64}~(a),
the disk is no longer a normal disk after the isotopy.  So, we can
perform another homotopy to make $F_k$ (after the first isotopy) a
normal surface.  This homotopy reduces $|F_k\cap\mathcal{T}^{(1)}|$ by
at least 2 as we push the disk in Figure~\ref{F64}~(a) across the
edge, which contradicts the assumption that $F_k$ has least weight.
Thus, every normal disk that contains such an arc $\delta$ is a
quadrilateral.  Since there are only two triangles in $\partial M$,
and since the edge $e$ lies inside $A$, there must be two arcs
$\delta_1$ and $\delta_2$ in $C$ that cut off two corners of the same
triangle (in the induced triangulation of $\partial M$).  By the
argument above, the two normal disks that contain $\delta_1$ and
$\delta_2$ respectively must be two quadrilaterals of different normal
disk types in the same tetrahedron.  Note that, during the isotopy as
in Figure~\ref{F63}, we push parts of $\partial F_k$ from
$N(B)\cap\partial M$ into the annulus $A$, and by unwrapping every
such spiral, we can assume that any two parallel normal disks in $F_k$
remain parallel after the isotopies.  We keep unwrapping the spirals
by isotopies as in Figure~\ref{F63}.  Either the weight of $F_k$ can
be reduced at a certain stage, or we can eventually push parts of
$\pi(D_k)\cap\partial M$ into the annulus $A$.  In particular, after
unwrapping the spirals enough times, we can assume that the $\delta_1$
and $\delta_2$ above lie in the cross disk.  Then, we can assume that
there is a pair of normal disks in the cross disk for each of the two
quadrilateral normal disk types that correspond to the $\delta_1$ and
$\delta_2$.  Since any two quadrilateral normal disks of different
types must intersect each other, those 4 quadrilaterals give rise to 4
planes in $\Tilde{F}_k$ intersecting each other (as in
Proposition~\ref{P41}), which contradicts the hypothesis that $F_k$
has the 4--plane property.  So, if $k$ is large enough, we can reduce
the weight of $F_k$ at a certain stage of the isotopy above.
Therefore, Figure~\ref{F61}~(b) cannot occur and claim \ref{claim1}
holds.
\end{proof}

\begin{figure}[ht!]
\cl{\relabelbox\small
\epsfxsize3.6in\epsfbox{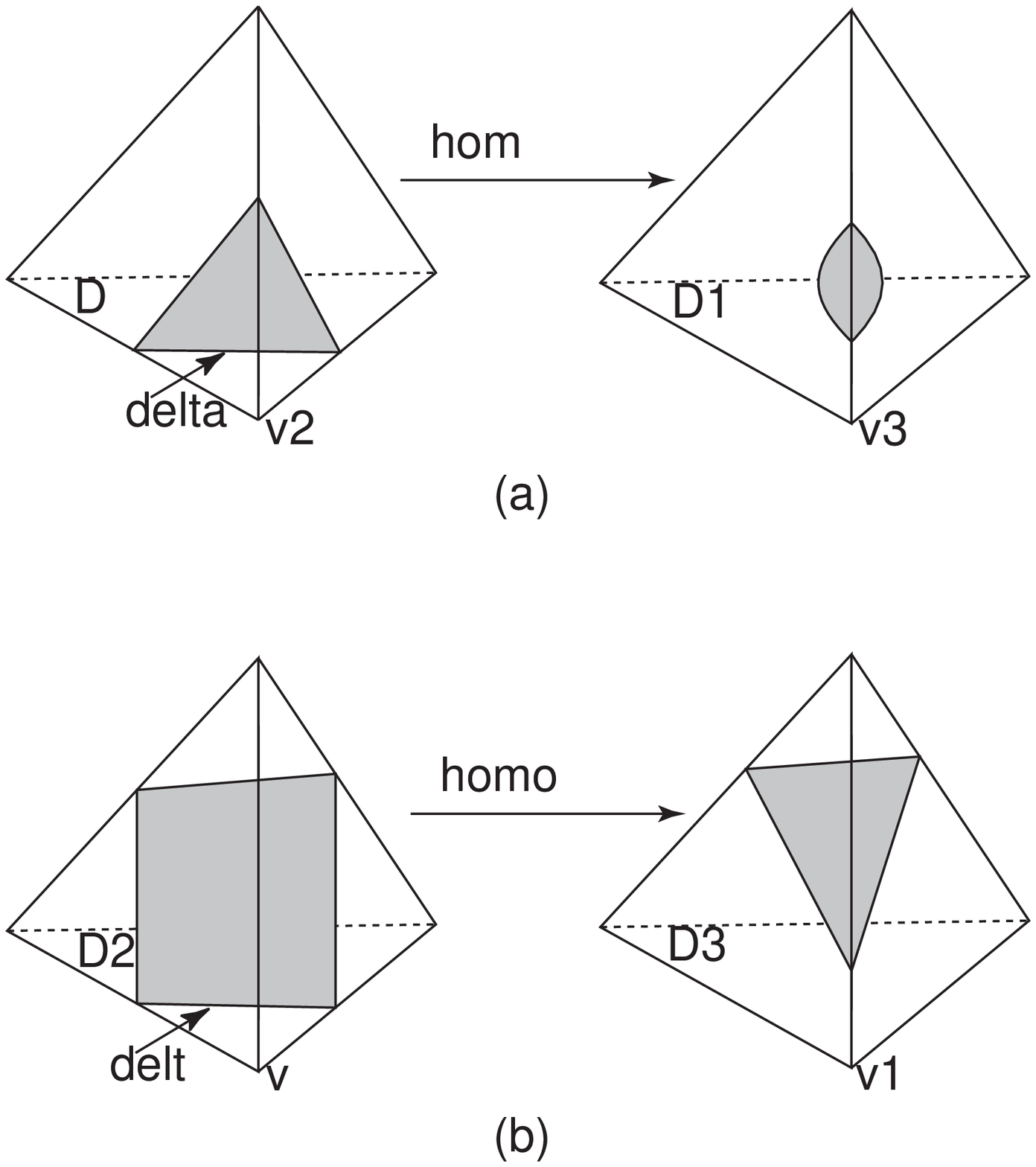}
\relabel {(a)}{(a)} 
\relabel {(b)}{(b)} 
\relabel {v}{$v$}
\relabel {v1}{$v$}
\relabel {v2}{$v$}
\relabel {v3}{$v$}
\relabel {hom}{homotopy} 
\relabel {homo}{homotopy} 
\relabel {delta}{$\delta$} 
\relabel {delt}{$\delta$} 
\relabel {D}{$\Delta$}
\relabel {D1}{$\Delta$}
\relabel {D2}{$\Delta$}
\relabel {D3}{$\Delta$}
\endrelabelbox}
\caption{}\label{F64}
\end{figure}

The branched surface that fully carries $\mu$ also fully carries a
compact surface, and by Lemma~\ref{L62}, the slope of $\partial\mu$ is
the same as the boundary slope of an incompressible and
$\partial$--incompressible surface. By Hatcher's theorem, there are
only finitely many possible slopes for $\partial\mu$. If the
lamination $\mu$ is constructed using cross disks from the sequence of
surfaces $\{F_k\}$, then the arcs in $\partial F_k$ must wind around
$\partial\mu$ many times (if $k$ is large). Therefore, by
Corollary~\ref{C32}, our construction of $\mu$ and Claim~\ref{claim1}
above, it is easy to see that there must be infinitely many slopes
that cannot be the boundary slopes of surfaces with the 4--plane
property, and Dehn fillings along these slopes yield 3--manifolds that
cannot admit any nonpositive cubings (see the proof of
Theorem~\ref{T3} at the end of this paper).  This can be viewed as a
weaker version of Theorems \ref{T2} and \ref{T3}.  To prove
Theorem~\ref{T2} to the full extent, which says that only finitely
many slopes can be the boundary slopes of surfaces with the 4--plane
property, we need to study the local pictures of the limit lamination
and surfaces with the 4--plane property in detail.

We first consider the case that $\mu$ is a compact orientable surface.
The proof for the case that $\mu$ contains a noncompact leaf is
similar.  Let $\mu\times I\subset M$ ($I=[0,1]$) be a small
neighborhood of $\mu$ in $M$, and $\tilde{\mu}\times I$ be a component
of the preimage of $\mu\times I$ in $\Tilde{M}$ with the induced fiber
structure.  Since $\mu$ is a compact embedded essential surface in
$M$, $\tilde{\mu}\times I$ can be considered as the universal cover of
$\mu\times I$, and we can assume $\pi_1(\mu)$ acts on $\tilde{M}$
fixing $\tilde{\mu}\times I$.  Suppose $k$ is large.  By our
construction of the lamination, there is always a large sub cross disk
of $D_k=D_k'\cup D_k''$ lying in $\tilde{\mu}\times I$.  To simplify
notation, we assume that $D_k\subset\tilde{\mu}\times I$; otherwise we
use a large sub cross disk of $D_k$ and the proof is the same.

Let $F_k'$ be the plane in $\Tilde{F}_k$ that contains $D_k'$,
$H'=F_k'\cap(\tilde{\mu}\times I)$, $H=\pi (H')$.  Since we can give
every component of $\partial F_k$ an orientation that agrees with the
induced orientation of $\partial F_k\cap(\mu\times I)$ in
claim~\ref{claim1}, we can assume the sign of every intersection point
of $\partial F_k\cap\partial S$ is always the same, where
$S=\mu\times\{t\}$ ($t\in I$).  Then, $H$ cannot be transverse to
every $I$--fiber of $\mu\times I$, because otherwise, by the argument
in the proof of Hatcher's theorem, $\partial F_k$ and $\partial S$
would have the same slope, which contradicts our assumptions.
Figure~\ref{F65} gives a local picture of $H$ where it is not
transverse to an $I$--fiber of $\mu\times I$.

In fact, it is not hard to see that, in some tetrahedron $T$, there
must be two different types of quadrilateral normal disks in $T\cap S$
and $T\cap F_k$ respectively.  Otherwise, by an argument in \cite{FO},
$H$ and $S$ lie in $N(B_T)$ and are transverse to the $I$--fibers of
$N(B_T)$, where $N(B_T)$ is a fibered neighborhood of an embedded
normal branched surface $B_T$.  Hence, by the arguments in the proof
of Hatcher's theorem, $F_k$ and $S$ have the same boundary slope
(although $F_k$ is not embedded), which contradicts our assumption.

Since all these surfaces are normal, after a small homotopy, we can
assume that each $I$--fiber of $\mu\times I$ either transversely
intersects $H$ or entirely lies in $H$, in which case the local
picture of this fiber is as shown in Figure~\ref{F65}, and we call
such fibers \emph{puncturing fibers}.  We can assume $\partial
F_k\cap(\partial\mu\times I)$ is a union of spirals, and by
claim~\ref{claim1}, the intersection points in $\partial
F_k\cap\partial S$ ($S=\mu\times\{t\}$) all have the same sign.  Then,
by our assumption on $H\cap(\mu\times I)$ and the argument above on
Hatcher's trick, any arc of $F_k\cap S$ with endpoints in $\partial M$
must pass through a puncturing fiber.  Since there is a large cross
disk wrapping around the compact surface $\mu$ many times, such a
puncturing fiber must puncture a cross disk, and we immediately get
three planes (in the universal cover) intersecting each other.
Moreover, any relatively short (compared with the size of the cross
disk) arc with endpoints in different components of $\mu\times\partial
I$ also punctures a cross disk.  Furthermore, if we can find two such
short arcs that are not far away from each other, then they puncture
the same cross disk.  If, in addition, the two planes that contain the
two short arcs intersect each other, we get a contradiction to the
4--plane property.  This is the basic idea of our proof.  After
perturbing $F_k$ a little, we can assume that $F_k\cap (\mu\times I)$
is transverse to the $I$--fibers of $\mu\times I$ except at puncturing
fibers and there are only finitely many puncturing fibers in $F_k\cap
(\mu\times I)$.

\begin{figure}[ht!]
\cl{\relabelbox\small
\epsfxsize3.6in\epsfbox{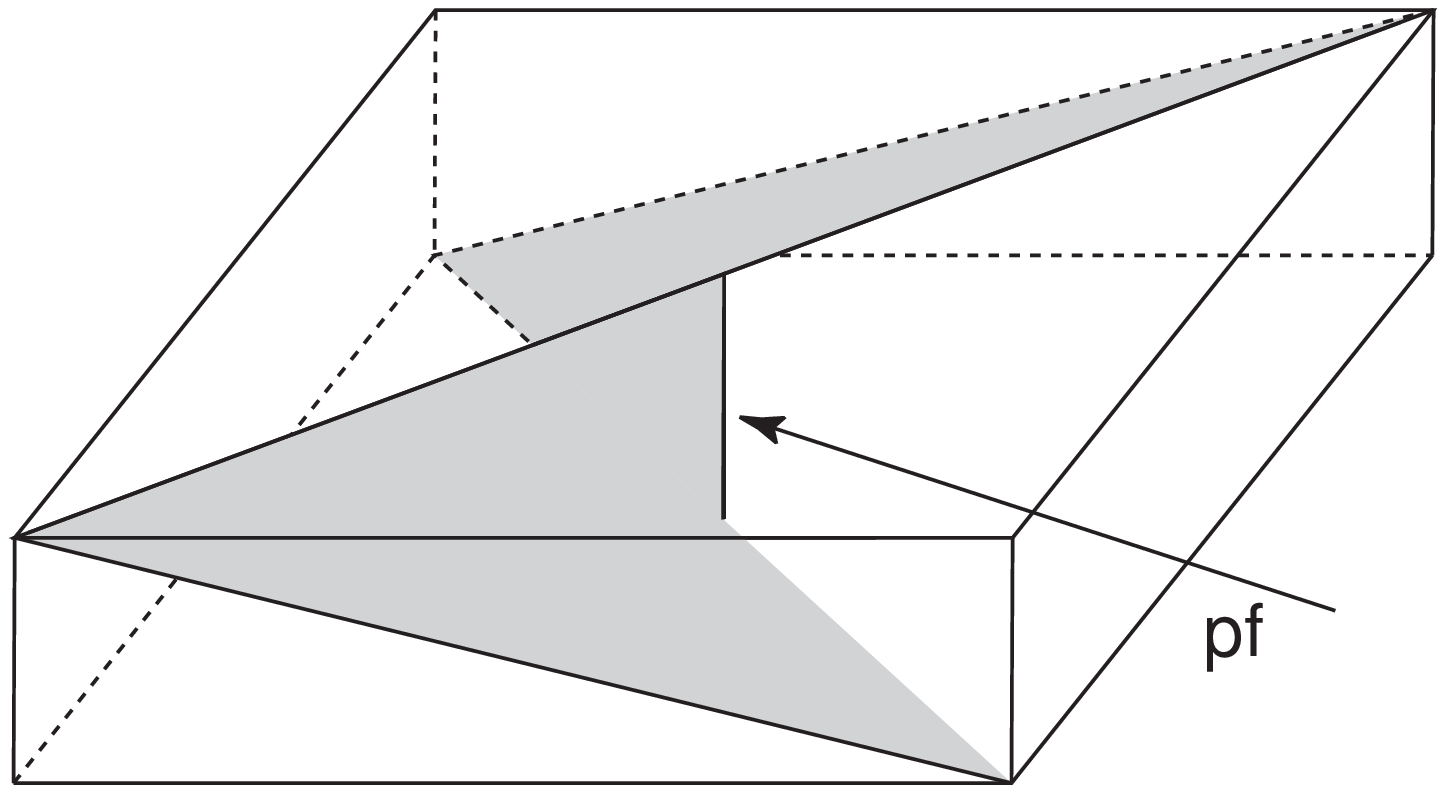}
\relabel {pf}{puncturing\,fiber} 
\endrelabelbox}
\caption{}\label{F65}
\end{figure}

The following observation, which summarizes the argument above, is
important to the remainder of the proof of Lemma~\ref{L63}.

\begin{observation}\label{principle}
Let $\alpha_i$ ($i=1,2$) be an arc in
$\tilde{F}_k\cap(\tilde{\mu}\times I)$ with two endpoints lying in
different components of $\tilde{\mu}\times\partial I$.  Suppose
$length(\alpha_1)$, $length(\alpha_2)$, and the distance between
$\alpha_1$ and $\alpha_2$ are bounded by a fixed number.  Then, if $k$
is large, $\alpha_1$ and $\alpha_2$ must puncture the same cross disk.
Let $F^{(i)}$ denote the plane in $\tilde{F}_k$ that contains
$\alpha_i$ ($i=1,2$).  If $F^{(1)}\cap F^{(2)}\ne\emptyset$ and
$F^{(1)}\ne F^{(2)}$, $F^{(1)}$, $F^{(2)}$ and the two planes
containing the two components of the cross disk are 4 planes in
$\tilde{F}_k$ intersecting each other, which contradicts the 4--plane
property.
\end{observation}

We denote the puncturing fibers of $F_k\cap(\mu\times I)$ by
$\gamma_1,\dots,\gamma_n$.  Let $q\co \mu\times I\to\mu$ be the
projection map.  Hence, $q(\gamma_1),\dots, q(\gamma_n)$ are points in
$\mu$.  Then, we can connect $q(\gamma_1),\dots, q(\gamma_n)$ buy
simple arcs $\alpha_1,\dots,\alpha_k$ to form a 1--complex $\Gamma$ in
which the $q(\gamma_i)$'s are the 0--cells and the $\alpha_i$'s are
the 1--cells.  Moreover, we can assume that $\mu-\Gamma$ is a union of
disks and annular neighborhoods of $\partial\mu$.  We denote the
closure of the annular components of $\mu-\Gamma$ by $A_1,\dots,A_m$,
where $m$ is the number boundary components of $\mu$.  Thus, for each
$i$, one component of $\partial A_i$ is a boundary circle of $\mu$ and
the other component of $\partial A_i$ lies in $\Gamma$.  We denote
$q^{-1}(\Gamma)$ and $q^{-1}(\alpha_i)$ by $\Gamma\times I$ and
$\alpha_i\times I$ respectively, and denote the preimage of
$\Gamma\times I$ in $\tilde{\mu}\times I$ by $\tilde{\Gamma}\times I$.

Let $S=\mu\times\{t\}$ and $\tilde{S}=\tilde{\mu}\times\{t\}$ ($t\in
I$).  We assume $S$ is transverse to $F_k$.  As before, any double arc
in $F_k\cap S$ or $\tilde{F}_k\cap\tilde{S}$ must pass through a
puncturing fiber.  Let $\beta$ be a subarc of a double arc of
$\tilde{F}_k\cap\tilde{S}$ in $\tilde{M}$ with one endpoint $Z_0$ in
$\partial\tilde{S}\cap\partial\tilde{F}_k$ and the other endpoint in a
puncturing fiber.  We can assume the interior of $\beta$ does not
intersect any puncturing fiber.  We denote the closure of the
component of $\tilde{\mu}\times I-\tilde{\Gamma}\times I$ containing
$Z_0$ by $\tilde{A}_1\times I$, and suppose $\pi(\tilde{A}_1\times I)$
is $A_1\times I$ in $\mu\times I$ defined above, where $\pi\co
\tilde{M}\to M$ is the covering map.

Let $D_{[k/2]}=D_{[k/2]}'\cup D_{[k/2]}''$ be a sub cross disk of
$D_k=D_k'\cup D_k''$ of size $[k/2]$ and with the same center.  By
choosing an appropriate $t\in I$, we can assume that $Z_0$ lies in
$D_{[k/2]}'$.  Moreover, after some isotopy on $F_k\cap(\mu\times I)$
(or choosing an appropriate $\Gamma$), we can assume that there is a
subarc of $\beta$, which we denote by $\beta_0$, properly embedded in
$\tilde{A}_1\times I$ such that $Z_0\in\partial\beta_0$ and
$\beta_0\subset D_{[k/2]}'$.  We denote the two lines in
$\partial\tilde{A}_1$ by $l_1$ and $l_2$, and suppose $Z_0\in
l_1\times I$.  Let $Z_1$ be the other endpoint of $\beta_0$.  Hence,
$Z_1\in l_2\times I\subset\tilde{\mu}\times I$.  By choosing $k$ to be
large, we can assume the length of the curve in
$\tilde{F}_k\cap(l_2\times I)$ that contains $Z_1$ is large, since
$Z_1$ lies in the sub cross disk $D_{[k/2]}$.

We can assume that the interior of $\beta$ is transverse to
$\tilde{\Gamma}\times I$.  Note that $int(\beta)$ does not intersect
any puncturing fiber in $\tilde{\Gamma}\times I$.  We denote the
points in $\beta\cap(\tilde{\Gamma}\times I)$ by $Z_1,\dots,Z_s$,
where $Z_1$ is as above and $Z_s=\partial\beta-Z_0$ lies in a
puncturing fiber.  These $Z_i$'s divide $\beta$ into $s$ subarcs
$\beta_0, \beta_1,\dots,\beta_{s-1}$, where $\partial\beta_i=Z_i\cup
Z_{i+1}$ and $\beta_0$ is as above.

We regard $\tilde{\Gamma}$ as a 1--complex in $\tilde{\mu}$ with
0--simplices corresponding to the puncturing fibers.  Let $\alpha$ be
any 1--simplex in $\tilde{\Gamma}$. So, $\alpha\times
I\subset\tilde{\Gamma}\times I$ is a vertical rectangle in
$\tilde{\mu}\times I$.  We call an arc in
$\tilde{F}_k\cap(\alpha\times I)$ a $\partial$--parallel arc if the
two endpoints of this arc lie in the same component of
$int(\alpha)\times\partial I$, where $int(\alpha)$ denotes the
interior of $\alpha$.  We can perform some normal homotopy on $F_k$ to
push all the $\partial$--parallel arcs out of $\pi(\alpha\times I)$,
where $\pi\co \tilde{M}\to M$ is the covering map, so that if two arcs
in $\tilde{F}_k\cap(\alpha\times I)$ do not intersect each other, then
after this homotopy, they do not intersect each other either.  Hence,
this normal homotopy preserve the 4--plane property.  Therefore, we
can assume that $\tilde{F}_k\cap(\alpha\times I)$ contains no
$\partial$--parallel arcs for any 1--simplex $\alpha$.  Moreover, if
$k$ is large, any $\partial$--parallel arc does not lie in the sub
cross disk $D_{[k/2]}$, and hence this homotopy does not affect the
previous assumptions on $\beta_0\subset D_{[k/2]}$.

Let $\zeta$ be an arc in $\tilde{F}_k\cap(\alpha\times I)$.  Since
there is no $\partial$--parallel arc in $\tilde{F}_k\cap(\alpha\times
I)$, either the two endpoints of $\zeta$ lie in different components
of $int(\alpha)\times\partial I$, or one endpoint of $\zeta$ lies in a
puncturing fiber $\gamma\subset\partial\alpha\times I$ in which case
we denote the two planes in $\tilde{F}_k$ containing $\zeta$ and
$\gamma$ by $F_{\zeta}$ and $F_{\gamma}$ respectively.  So, if the
second case happens, either $F_{\zeta}\ne F_{\gamma}$ and
$F_{\zeta}\cap F_{\gamma}\ne\emptyset$ or $F_{\zeta}=F_{\gamma}$.  We
call $\zeta$ a \emph{puncturing arc} if either the two endpoints of
$\zeta$ lie in different components of $int(\alpha)\times\partial I$,
or $F_{\zeta}=F_{\gamma}$.  Thus, if $\zeta$ is a puncturing arc,
there must be a relatively short arc in $F_{\zeta}$ containing $\zeta$
and with two endpoints in different components of
$\tilde{\mu}\times\partial I$.  The role of a puncturing arc is the
same as the role of a puncturing fiber, see
Observation~\ref{principle}.  Moreover, if $\zeta$ is not a puncturing
arc, then one of the two puncturing fibers in $\partial\alpha\times I$
intersects the plane $F_{\zeta}$ nontrivially.

\begin{claim}~\label{pa}
Each $Z_i\in\beta$ ($1\le i\le s-1$) lies in a puncturing arc.
\end{claim}
\begin{proof}[Proof of Claim~\ref{pa}]
We first show that $Z_{s-1}$ lies in a puncturing arc, and then we
inductively prove it for each $Z_i$.  Suppose $Z_{s-1}$ lies in
$\alpha\times I$, where $\alpha$ is a 1--simplex of $\tilde{\Gamma}$,
and we denote the arc in $\tilde{F}_k\cap(\alpha\times I)$ containing
$Z_{s-1}$ by $\zeta_{s-1}$.

Let $A\times I\subset\tilde{\mu}\times I$ be the closure of the
component of $\tilde{\mu}\times I-\tilde{\Gamma}\times I$ that
contains $\beta_{s-1}$.  So, $\beta_{s-1}$ is an properly embedded arc
in $A\times I$ with $\partial\beta_{s-1}\subset\partial A\times I$.
Since $\mu-\Gamma$ consists of disks and annular neighborhoods of
circles in $\partial\mu$, either $A$ is a disk, or $A$ is of the form
$[0,1]\times\mathbb{R}$ which can be considered as the universal cover
of an annular neighborhood of a boundary circle of $\mu$.

If $A$ is a disk, since $\beta_{s-1}$ is properly embedded in $A\times
I$, after some isotopy on $F_k$ if necessary, we can assume
$length(\beta_{s-1})$ is bounded by a number that depends only on $A$.
If $\zeta_{s-1}$ is not a puncturing arc in $\alpha\times I$, since we
have assumed that $\tilde{F}_k\cap(\alpha\times I)$ does not contain
any $\partial$--parallel arc, $\zeta_{s-1}$ must intersect at least
one of the two puncturing fibers in $\partial\alpha\times I$.  Let
$\gamma$ be a puncturing fiber in $\partial\alpha\times I$ that
intersects $\zeta_{s-1}$.  By our definition of puncturing arc, the
two planes in $\tilde{F}_k$ containing $\gamma$ and $\beta_{s-1}$
intersect each other.  By our construction above,
$Z_s\in\partial\beta_{s-1}$ lies in another puncturing fiber, say
$\gamma_s$.  Since $A$ is a disk, the distance between the two
puncturing fibers $\gamma$ and $\gamma_s$ is bounded by the diameter
of the disk $A$.  So, $\gamma$ and $\gamma_s$ puncture the same cross
disk.  Moreover, the plane containing $\gamma_s$, $\beta_{s-1}$ and
$\zeta_{s-1}$ and the plane containing $\gamma$ intersect each other,
which contradicts the 4--plane property as in
Observation~\ref{principle}.

So, we only need to consider the case that $A$ is of the form
$[0,1]\times\mathbb{R}$.  We denote the two boundary lines of $A$ by
$l_i=\{i\}\times\mathbb{R}$ ($i=0,1$).  Suppose
$l_0\subset\partial\tilde{\mu}$.  Hence, $\pi(l_0)$ is a boundary
circle of $\mu$, where $\pi\co \tilde{M}\to M$ is the covering map.
This circle $\pi(l_0)$ represents a nontrivial element
$g\in\pi_1(\mu)\subset\pi_1(M)$, and $g$ acts on $\tilde{M}$ fixing
$A\times I$.  By our construction, unless $s=1$, $\beta_{s-1}$ is an
arc properly embedded in $A\times\{t\}$ with both endpoints in
$l_1\times\{t\}$.

If the length of the subarc of $l_1\times\{t\}$ between $Z_s$ and
$Z_{s-1}$ is large, then $g(Z_s)$ (or $g^{-1}(Z_s)$) lies between
$Z_s$ and $Z_{s-1}$ in $l_1\times\{t\}$, and hence $\beta_{s-1}$
nontrivially intersects $g(\beta_{s-1})$ (or $g^{-1}(\beta_{s-1})$).
Note that $Z_s$ and $g(Z_s)$ lie in puncturing fibers.  Since
$g\in\pi_1(\mu)$ is fixed, the distance between the two puncturing
fibers containing $Z_s$ and $g(Z_s)$ is relatively small (compared
with $k$), and hence they puncture the same cross disk.  Moreover,
since $\beta_{s-1}$ nontrivially intersects $g(\beta_{s-1})$ (or
$g^{-1}(\beta_{s-1})$), the two planes containing $Z_s$ and $g(Z_s)$
(or $g^{-1}(Z_s)$) intersect each other, which contradicts the
4--plane property as in Observation~\ref{principle}.

Thus, we can assume length of the subarc of $l_1\times\{t\}$ between
$Z_s$ and $Z_{s-1}$ is relatively small.  Hence, the distance between
$Z_s$ and the two puncturing fibers $\partial\alpha\times I$ is
relatively small, where $\alpha$ is a 1--simplex of $\tilde{\Gamma}$
and $Z_{s-1}\in\alpha\times I$.  As in Observation~\ref{principle},
the plane in $\tilde{F}_k$ containing $Z_s\cup\beta_{s-1}$ and a plane
containing a puncturing fiber in $\partial\alpha\times I$ cannot
intersect each other.  So, as in the case that $A$ is a disk, the arc
$\zeta_{s-1}$, which is the arc in $\tilde{F}_k\cap(\alpha\times I)$
containing $Z_{s-1}$, must be a puncturing arc in $\alpha\times I$.
Therefore, in any case, $Z_{s-1}$ lies in a puncturing arc.

Then, we apply the argument above to $\beta_{s-2}$ ($s>2$).  Now,
$\partial\beta_{s-2}=Z_{s-1}\cup Z_{s-2}$.  Note that in this case
$Z_{s-1}$ lies in a puncturing arc ($Z_s$ lies in a puncturing fiber
in the case above), but this does not make any difference when using
Observation~\ref{principle}.  Hence, $Z_{s-2}$ lies in a puncturing
arc, and inductively, each $Z_i\in\beta$ ($1\le i\le s-1$) lies in a
puncturing arc.
\end{proof}

By our assumption before, $Z_1$ lies in the sub cross disk
$D_{[k/2]}$.  We can choose $k$ large enough so that there is no short
arc containing $Z_1$ and with endpoints in different components of
$\tilde{\mu}\times\partial I$. This contradicts Claim~\ref{pa}.  Thus,
Lemma~\ref{L63} holds in the case that $\mu$ is a compact orientable
surface.

If $\mu$ is a compact nonorientable surface, since $M$ is orientable,
we can apply Hatcher's trick to the horizontal boundary of a twisted
$I$--bundle over $\mu$, and the proof is the same.

Suppose $\mu$ contains a noncompact leaf.  Let $B$ be a branched
surface fully carries $\mu$, $L$ be the branch locus of $B$, and $p\co
N(B)\to B$ be the map that collapses every $I$--fiber of $N(B)$ to a
point.  We can assume $\partial B$ is a union of circles in $\partial
M$.  By previous arguments, any such branched surface always fully
carries a compact surface with the same boundary slope as
$\partial\mu$.  Let $S$ be a compact surface fully carried by $B$.  By
Claim~\ref{claim1} and Hatcher's trick, as in the case that $\mu$ is a
compact surface, $F_k$ is not transverse to the $I$--fibers of $N(B)$
along any arc of $F_k\cap S$.  As before, in some tetrahedron $T$,
there must be two different types of quadrilateral normal disks in
$T\cap S$ and $T\cap F_k$ respectively.  Thus, after a small homotopy,
we can assume that each $I$--fiber of $N(B)$ either transversely
intersects $F_k$ or entirely lies in $F_k$, in which case the local
picture of this fiber is as shown in Figure~\ref{F65} and we also call
such fibers puncturing fibers.  We can assume there are only finitely
many puncturing fibers for each $F_k$.

$N(B)$ can be viewed as the gluing of a collection of $I$--bundles
over compact surfaces along $p^{-1}(L)$.  Now, we use the puncturing
fibers to decompose $N(B)$ into a similar structure.  We say
$\alpha\times I\subset N(B)$ is \emph{vertical} if $\{p\}\times I$ is
a subarc of an $I$--fiber of $N(B)$ for each $p\in\alpha$ and
$(\alpha\times I)\cap\mu=\alpha\times C$ for some closed set $C\subset
I$.  We start with the puncturing fibers of $F_k$.  Since every leaf
of $\mu$ is dense, we can add finitely many vertical rectangles
$\alpha_i\times I\subset N(B)$ ($i=1,\dots,n$) such that
$\partial\alpha_i\times I$ is a pair of subarcs of puncturing fibers
for each $i$ and $\mu-\cup_{i=1}^n\alpha_i\times I$ consists of disks
and annular neighborhoods of circles in $\partial\mu$.  Moreover, we
can assume that there is a union of products $d_i\times I$
($i=1,\dots,s$) and $A_i\times I$ ($i=1,\dots,t$) that are glued along
$\cup_{i=1}^n\alpha_i\times I$, such that:
\begin{enumerate}
\item each $d_i$ is a disk and $\partial d_i\times I$ lies in
$\cup_{i=1}^n\alpha_i\times I$ for each $i$;
\item each $A_i$ is an annulus, one component of $\partial A_i\times
I$ lies in $\partial M$ and the other component of $\partial A_i\times
I$ lies in $\cup_{i=1}^n\alpha_i\times I$;
\item $\{p\}\times I$ is a subarc of an $I$--fiber of $N(B)$ for each
$p$ in $d_i$ or $A_i$;
\item $\mu$ lies in the union of these products $\alpha_i\times I$'s,
$d_i\times I$'s and $A_i\times I$'s;
\item $\mu\cap (d_i\times I)=d_i\times C_i$ and $\mu\cap (A_i\times
I)=A_i\times C_i'$, where $C_i$ and $C_i'$ are closed sets in
$int(I)$.
\end{enumerate}
Furthermore, we can assume the diameter of $d_i$ and the length of
$\partial A_i$ are bounded by a number independent of the puncturing
fibers, since $\mu$ is fixed.  In fact, after a small perturbation, we
can view the union of these products $\alpha_i\times I$'s, $d_i\times
I$'s and $A_i\times I$'s as a fibered neighborhood $N(B')$ of another
branched surface $B'$ that also fully carries $\mu$.  We can also view
$\cup_{i=1}^n\alpha_i\times I$ as $p^{-1}(L')$, where $L'$ is the
branch locus of $B'$ and $p\co N(B')\to B'$ is the map collapsing
every $I$--fiber to a point.

The new branched surface $B'$ also fully carries a compact surface,
say $S$.  We can suppose $S$ lies in $N(B')$ and $S$ does not
intersect $d_i\times\partial I$ or $A_i\times\partial I$. By
claim~\ref{claim1}, we can assume each point in $\partial
S\cap\partial F_k$ has the same sign.  Since the size of each $d_i$ is
bounded by a number independent of the puncturing fibers, if $k$ is
large, there is a cross disk from $F_k$ cutting through $d_i\times I$
for each $i$. Moreover, we can choose an appropriate surface $S$ so
that at least one point of $\partial S\cap\partial F_k$ belongs to a
sub cross disk $D_{[k/2]}$ as before.  By our construction, $\partial
d_i\times I$ and $\partial A_i\times I$ contain subarcs of puncturing
fibers which puncture a cross disk. After some homotopy as in the case
that $\mu$ is a compact surface, we can also assume that any
$\partial$--parallel arc of $F_k$ in $\partial d_i\times I$ or
$\partial A_i\times I$ does not intersect $S$.  Then, we can define
puncturing arcs using the $d_i\times I$'s and $A_i\times I$'s similar
to the case that $\mu$ is a compact surface, and the proof is the
same.
\end{proof}

Theorem~\ref{T2}, which is a generalization of Hatcher's theorem, now
follows easily from Corollary~\ref{C32} and Lemmas \ref{L61} and
\ref{L63}.

\begin{theorem2}
Let $M$ be an orientable and irreducible 3--manifold whose boundary is
an incompressible torus, and let $\mathcal{H}$ be the set of injective
surfaces that are embedded along their boundaries and satisfy the
4--plane property.  Suppose that $M$ does not contain any
nonperipheral closed (embedded) incompressible surfaces.  Then the
surfaces in $\mathcal{H}$ can realize only finitely many slopes.
\end{theorem2}
\begin{proof}
Suppose that the surfaces can realize infinitely many slopes.  Let $\{
F_n\}$ be a sequence of surfaces in $\mathcal{H}$ no two of which have
the same boundary slopes.  Since they have different boundary slopes,
by Corollary~\ref{C32}, the surfaces in $\{ F_n\}$ cannot be fully
carried by finitely many immersed branched surfaces.  Then, by the
argument in section~\ref{S4}, there exist a sequence of cross disks
from $\{ F_n\}$ that gives rise to an essential lamination.  However,
Lemma~\ref{L63} imply that the sequence $\{ F_n\}$ contains a
subsequence of surfaces with the same boundary slope, which
contradicts our assumption that the surfaces in $\{ F_n\}$ all have
different boundary slopes.
\end{proof}

As an application of Theorem~\ref{T2}, we prove Theorem~\ref{T3},
which gives the first nontrivial examples of 3--manifolds that do not
admit any nonpositive cubings.  Before we proceed, we prove the
following lemma.

\begin{lemma}\label{L64}
Let $M$ be a closed and irreducible 3--manifold, $S$ be a closed least
weight surface in $M$ with the 4--plane property, and $C$ be a
homotopically nontrivial simple closed curve that intersects $S$
nontrivially.  Then $S-C$ is a surface with the 4--plane property in
$M-C$.
\end{lemma}
\begin{proof}
Let $\Tilde{M}$ be the universal cover of $M$ and $\Tilde{C}$ be the
preimage of $C$ in $\Tilde{M}$.  So, $\Tilde{M}-\Tilde{C}$ is a cover
of $M-C$.  Let $\tilde{S}$ be the preimage of $S$ in $\tilde{M}$.
Then, $\Tilde{S}-\Tilde{C}$ is a collection of embedded surfaces in
$\tilde{M}-\tilde{C}$.  Since $S$ has the 4--plane property, among any
$4$ embedded surfaces in $\Tilde{S}-\Tilde{C}$, there is a disjoint
pair.  Moreover, as each surface in $\Tilde{S}-\Tilde{C}$ is embedded,
among any $4$ planes in the preimage of $\Tilde{S}-\Tilde{C}$ in the
universal cover of $\Tilde{M}-\Tilde{C}$ (ie the universal cover of
$M-C$), there is a disjoint pair.  Therefore, $S-C$ satisfies the
$4$--plane property in $M-C$.
\end{proof}

\begin{theorem3}
Let $M$ be an orientable and irreducible 3--manifold whose boundary is
an incompressible torus.  Suppose that $M$ does not contain any closed
nonperipheral (embedded) incompressible surfaces.  Then only finitely
many Dehn fillings on $M$ can yield 3--manifolds that admit
nonpositive cubings.
\end{theorem3}
\begin{proof}
Let $M(s)$ be the closed 3--manifold after doing Dehn filling along
slope $s$, and $C_s$ be the core of the solid torus glued to $M$
during the Dehn filling.  Then, except for finitely many slopes, $C_s$
is a homotopically nontrivial curve in $M(s)$.  Suppose that $M(s)$
admits a nonpositive cubing.  For each cube in the cubing, there are
$3$ disks parallel to the square faces and that intersect the edges of
the cube in their mid-points.  These mid-disks from all the cubes in
the cubing match up and yield a union of immersed surfaces, which we
denote by $\mathcal{S}$.  The complement of $\mathcal{S}$ is a union
of 3--balls.  Aitchison and Rubinstein have shown that these surfaces
(and their double covers in $M(s)$ if they are one-sided) satisfy the
$4$--plane property \cite{AR}.  Since $C_s$ is nontrivial and the
complement of $\mathcal{S}$ is a union of 3--balls, $C_s$ must
nontrivially intersect at least one immersed surface in $\mathcal{S}$.
Let $N(C_s)$ be a small tubular neighborhood of $C_s$.  Note that
$\mathcal{S}-int(N(C_s))$ may not be $\partial$--injective in $M$ and
we need to perform some homotopy on the surfaces in $\mathcal{S}$.
Similar to the case of embedded incompressible surfaces \cite{Th}, we
push the (immersed) $\partial$--compressing disk across $N(C_s)$ and
get less intersection circles.  Since $\mathcal{S}$ is immersed, this
homotopy changes the intersection patterns of $\mathcal{S}$ in $M(s)$,
but by choosing innermost $\partial$--compressing disks if necessary,
we can require that the disjoint planes in the preimage of
$\mathcal{S}$ in the universal cover of $M(s)$ remain disjoint after
this homotopy, and hence this homotopy preserves the 4--plane property
of $\mathcal{S}$.

The nonpositive cubing gives $M(s)$ a singular nonpositive metric and
$\mathcal{S}$ consists of totally geodesic surfaces in this singular
metric \cite{AR}.  The geodesic that represents $C_s$ must intersect
some (totally geodesic) surface in $\mathcal{S}$.  Since the singular
metric is nonpositive, after the homotopy above, $C_s$ must still
intersect some immersed surface in $\mathcal{S}$.  Hence, by
Lemma~\ref{L64}, there is an injective surface in $M$ that satisfies
the $4$--plane property and has boundary slope $s$.  By
Theorem~\ref{T2}, there are only finitely many such slopes.
Therefore, the theorem holds.
\end{proof}

\end{document}